\documentclass[11pt]{article}
\usepackage{amssymb,amsmath,latexsym}
\usepackage{pdfsync}

\allowdisplaybreaks

\begin{document}

\begin{doublespace}

\def\1{{\bf 1}}
\def\ind{{\bf 1}}
\def\nn{\nonumber}

\def\sA {{\cal A}} \def\sB {{\cal B}} \def\sC {{\cal C}}
\def\sD {{\cal D}} \def\sE {{\cal E}} \def\sF {{\cal F}}
\def\sG {{\cal G}} \def\sH {{\cal H}} \def\sI {{\cal I}}
\def\sJ {{\cal J}} \def\sK {{\cal K}} \def\sL {{\cal L}}
\def\sM {{\cal M}} \def\sN {{\cal N}} \def\sO {{\cal O}}
\def\sP {{\cal P}} \def\sQ {{\cal Q}} \def\sR {{\cal R}}
\def\sS {{\cal S}} \def\sT {{\cal T}} \def\sU {{\cal U}}
\def\sV {{\cal V}} \def\sW {{\cal W}} \def\sX {{\cal X}}
\def\sY {{\cal Y}} \def\sZ {{\cal Z}}

\def\bA {{\mathbb A}} \def\bB {{\mathbb B}} \def\bC {{\mathbb C}}
\def\bD {{\mathbb D}} \def\bE {{\mathbb E}} \def\bF {{\mathbb F}}
\def\bG {{\mathbb G}} \def\bH {{\mathbb H}} \def\bI {{\mathbb I}}
\def\bJ {{\mathbb J}} \def\bK {{\mathbb K}} \def\bL {{\mathbb L}}
\def\bM {{\mathbb M}} \def\bN {{\mathbb N}} \def\bO {{\mathbb O}}
\def\bP {{\mathbb P}} \def\bQ {{\mathbb Q}} \def\bR {{\mathbb R}}
\def\bS {{\mathbb S}} \def\bT {{\mathbb T}} \def\bU {{\mathbb U}}
\def\bV {{\mathbb V}} \def\bW {{\mathbb W}} \def\bX {{\mathbb X}}
\def\bY {{\mathbb Y}} \def\bZ {{\mathbb Z}}
\def\R {{\mathbb R}} \def\RR {{\mathbb R}}
\def\n{{\bf n}}

\newcommand{\expr}[1]{\left( #1 \right)}
\newcommand{\cl}[1]{\overline{#1}}
\newtheorem{thm}{Theorem}[section]
\newtheorem{lemma}[thm]{Lemma}
\newtheorem{defn}[thm]{Definition}
\newtheorem{prop}[thm]{Proposition}
\newtheorem{corollary}[thm]{Corollary}
\newtheorem{remark}[thm]{Remark}
\newtheorem{example}[thm]{Example}
\numberwithin{equation}{section}
\def\ee{\varepsilon}
\def\qed{{\hfill $\Box$ \bigskip}}
\def\NN{{\cal N}}
\def\AA{{\cal A}}
\def\MM{{\cal M}}
\def\BB{{\cal B}}
\def\CC{{\cal C}}
\def\LL{{\cal L}}
\def\DD{{\cal D}}
\def\FF{{\cal F}}
\def\EE{{\cal E}}
\def\QQ{{\cal Q}}
\def\RR{{\mathbb R}}
\def\R{{\mathbb R}}
\def\L{{\bf L}}
\def\K{{\bf K}}
\def\S{{\bf S}}
\def\A{{\bf A}}
\def\E{{\mathbb E}}
\def\F{{\bf F}}
\def\P{{\mathbb P}}
\def\N{{\mathbb N}}
\def\eps{\varepsilon}
\def\wh{\widehat}
\def\wt{\widetilde}
\def\pf{\noindent{\bf Proof.} }
\def\pff{\noindent{\bf Proof} }
\def\beq{\begin{equation}}
\def\eeq{\end{equation}}
\def\bee{\begin{equation}}
\def\eee{\end{equation}}
\def\osc{\mathrm{Osc}}

\title{\Large \bf
Global uniform boundary Harnack principle with explicit decay rate and its
application}
\author{{\bf Panki Kim}\thanks{This work was supported by Basic Science Research Program through the National Research Foundation of Korea(NRF) Grant funded by the Korea government(MEST)
(2012-0000940).} \quad {\bf Renming Song}\thanks{Research supported in part by a grant from the Simons
Foundation (208236).}
\quad and \quad {\bf Zoran Vondra\v{c}ek}\thanks{Supported in part by the MZOS
grant 037-0372790-2801.}  }

\date{ }
\maketitle

\begin{abstract}
In this paper, we consider a large class of
subordinate Brownian motions $X$ via subordinators with Laplace
exponents which are complete Bernstein functions satisfying
some mild scaling conditions at zero and at infinity.
We first discuss how such conditions govern the
behavior of the subordinator and the corresponding subordinate
Brownian motion for both large and small time and space.
Then we establish a global uniform
boundary Harnack principle in (unbounded) open sets for the subordinate Brownian motion.
When the open set satisfies the interior and exterior ball conditions with radius $R >0$,
we get a global uniform boundary Harnack principle with explicit decay rate.
Our boundary Harnack principle is global in the sense that it holds
for all $R>0$ and the comparison  constant  does not depend on  $R$, and it is uniform
in the sense that it holds for all balls with radii $r \le R$ and
the comparison  constant  depends neither on $D$ nor on $r$.
As an application, we give sharp two-sided estimates for the
transition densities and Green functions of such subordinate Brownian motions
in the half-space.
\end{abstract}



\section{Introduction}
The study of potential theory of discontinuous L\'evy processes in $\R^d$ revolves around several
fundamental questions such as sharp heat kernel and Green function estimates,
exit time estimates and Poisson kernel estimates, Harnack and boundary Harnack principles for
non-negative harmonic functions.
One can roughly divide these studies in two categories:
those on a bounded set and
those on an unbounded set.
For the former, it is the local behavior of the process that matters,
while for the latter both local and global
behaviors are important.
The processes investigated in these studies are usually described in two ways:
either the process is given explicitly through its characteristic
exponent (such as the case of a symmetric stable process, a relativistically stable process,
sum of two independent stable processes, etc.), or some conditions on the
characteristic exponent are given. In the situation when one is interested
in the potential theory on bounded sets, conditions imposed on the characteristic
exponent govern the small time -- small space (i.e., local) behavior of the process.
Let us be more precise and describe in some detail one such condition and
some of the results in the literature.

Let $S=(S_t)_{t\ge 0}$ be a subordinator  (that is, an
increasing L\'evy process satisfying $S_0=0$) with Laplace exponent $\phi$,
and let $W=(W_t)_{t\ge 0}$ be a Brownian motion in $\R^d$, $d\ge 1$,  independent of $S$ with
$$
\E_x\left[e^{i\xi(W_t-W_0)}\right]=e^{-t{|\xi|^2}}\ , \quad \xi\in \R^d, t>0.
$$
The process $X=(X_t)_{t\ge 0}$ defined by $X_t:=W(S_t)$ is called a
subordinate Brownian motion. It is a rotationally invariant L\'evy process in $\R^d$
with characteristic exponent $\phi(|\xi|^2)$.
The function $\phi$ is a Bernstein function.
Let us introduce the following upper and lower scaling conditions:

\medskip
\noindent
{\bf (H1):}
There exist constants $0<\delta_1\le \delta_2 <1$ and $a_1, a_2>0$  such that
\begin{equation}\label{e:H1}
a_1\left(\frac{R}{r}\right)^{\delta_1}\le \frac{\phi(R)}{\phi(r)}\le a_2\left(\frac{R}{r}\right)^{\delta_2},
\quad 1\le r\le R.
\end{equation}
It follows from the definitions in \cite[pp. 65 and 68]{BGT} and \cite[Proposition 2.2.1]{BGT} that
\eqref{e:H1} is equivalent to saying that $\phi$ is in the class $OR$ of $O$-regularly varying functions at $\infty$ with
Matuszewka indices contained in $(0, 1)$. The advantage of the formulation above is that we can provide
more direct proofs for some of the results below.
\eqref{e:H1} is a condition on the asymptotic behavior of $\phi$ at infinity and
it governs the behavior of the subordinator $S$
for small time and small space, which, in turn, implies the small time -- small space
behavior of the corresponding subordinate Brownian motion $X$. Very recently
it has been shown in \cite{KSV7} (see also \cite{KM2}) that if {\bf (H1)} holds and $\phi$ is a
complete Bernstein function, then the uniform boundary Harnack principle is true and
various exit time and Poisson kernel estimates hold.
Further, sharp two-sided Green function estimates for bounded $C^{1,1}$ open sets are given in \cite{KM2}.
The statements of these results usually take the following form:
For some $R>0$, there exists a constant $c=c(R)>0$ (also depending on the process $X$)
such that some quantities involving $r\in (0, R)$ can be estimated by expressions involving the constant $c$.
The point is that although the constant $c$ is uniform for small $r\in (0, R)$, it does depend on $R$,
meaning that the result is local. It would be of interest to obtain
a global and uniform version of such results,
namely with the constant depending neither on $R$ nor on the open set itself.
This would facilitate the study of potential theory on unbounded sets.
In order to accomplish this goal, it is
clear that the assumption {\bf (H1)} (or some similar condition) will not suffice,
and that one needs additional assumptions that govern the behavior of the process for
large time and large space.

In some recent papers (see \cite{CKS4, CKS7, CT})
potential-theoretic properties of stable and relativistically stable processes are
studied in unbounded sets such as the half-space, half-space-like $C^{1,1}$ open
sets and exterior $C^{1,1}$ open sets. Note that these processes are given
explicitly by its characteristic exponent. In the current paper we would like to
impose a condition similar to {\bf (H1)} that governs the large time -- large space
behavior of the process and
obtain global uniform potential-theoretic results.
Thus, in addition to {\bf (H1)}, we will also assume

\medskip
\noindent
{\bf (H2):}
There exist constants $0<\delta_3\le \delta_4 <1$ and $a_3, a_4>0$  such that
\begin{equation}\label{e:H2}
a_3\left(\frac{R}{r}\right)^{\delta_3}\le \frac{\phi(R)}{\phi(r)}\le a_4\left(\frac{R}{r}\right)^{\delta_4},
\quad r\le R\le 1.
\end{equation}
Similarly, \eqref{e:H2} is equivalent to saying that $\phi$ is in the class of
$O$-regularly varying functions
at 0 with Matuszewka indices contained in $(0, 1)$.
\eqref{e:H2} is a condition about the asymptotic behavior of $\phi$ at zero
and it governs the behavior of the
subordinator $S$ and the corresponding subordinate Brownian motion $X$ for large time and large space.
Also note that under {\bf (H2)}, $X$ is transient if $d\ge 2$.

Throughout the paper we will assume that $\phi$ is a complete Bernstein function satisfying
{\bf (H1)} and/or {\bf (H2)}, and $X_t=W(S_t)$ will be the
corresponding subordinate Brownian motion. First we study consequences of
scaling conditions on the subordinator $S$, its L\'evy density and potential density.
This is done in Section 2 of the paper. In Section 3 we proceed to
properties of the subordinate Brownian motion $X$. The first main result
is about estimates of the L\'evy density and the Green function of $X$ for the whole space
given in Theorem \ref{t:J-G}. These estimates allow us to repeat arguments from
\cite{KSV3, KSV7} and obtain global uniform estimates of the exit times and
the Poisson kernel, as well as global uniform Harnack and boundary Harnack principles.
The latter will play a crucial role in this paper.

In Section 4 we prove the main result of the paper -- the global uniform boundary Harnack principle
with explicit decay rate in open sets satisfying both interior and exterior ball conditions
(see Theorem \ref{L:2}(b)). The key technical contribution is Proposition \ref{l:main} which has
appeared in similar forms in several recent papers. The main novelty of the current version
is that the estimate gets better as the radius grows larger.
The quite technical part of the proof of this proposition is given in two auxiliary lemmas.

Theorem \ref{L:2} is used in Section 5 to obtain sharp two-sided heat kernel and Green function
estimates for the process $X$ killed upon exiting the half-space
$\bH=\{x=(x_1,\dots, x_{d-1},x_d)\in \R^d:\, x_d>0\}$. To the best of
our knowledge, this is the first time the heat kernel estimates are obtained
in an unbounded set for a process which is not given by an explicit characteristic exponent.

The results of this paper, especially those of Sections 3 and 4, are used in
the subsequent paper \cite{KSV9} to prove the boundary Harnack principle at infinity.
This was the main motivation for the investigations in the current paper.

Using the tables at the end of \cite{SSV}, one can construct a lot of explicit
examples of complete Bernstein functions satisfying both {\bf (H1)} and {\bf (H2)}.
Here are a few of them:
\begin{description}
\item{(1)} $\phi(\lambda)=\lambda^\alpha + \lambda^\beta$, $0<\alpha<\beta<1$;
\item{(2)} $\phi(\lambda)=(\lambda+\lambda^\alpha)^\beta$, $\alpha, \beta\in (0, 1)$;
\item{(3)} $\phi(\lambda)=\lambda^\alpha(\log(1+\lambda))^\beta$, $\alpha\in (0, 1)$,
$\beta\in (0, 1-\alpha)$;
\item{(4)} $\phi(\lambda)=\lambda^\alpha(\log(1+\lambda))^{-\beta}$, $\alpha\in (0, 1)$,
$\beta\in (0, \alpha)$;
\item{(5)} $\phi(\lambda)=(\log(\cosh(\sqrt{\lambda})))^\alpha$, $\alpha\in (0, 1)$;
\item{(6)} $\phi(\lambda)=(\log(\sinh(\sqrt{\lambda}))-\log\sqrt{\lambda})^\alpha$, $\alpha\in (0, 1)$.
\end{description}
We remark here that relativistic stable processes do not satisfy {\bf (H2)}, so
the present paper does not cover this interesting case. We plan to address this
important case in the near future.

Throughout this paper, $d\geq 1$ and  the constants
$C_1$, $a_i$ and $\delta_i$, $i=1, \dots, 4$,  will be fixed.
  We use $c_1, c_2, \dots$ to
denote generic constants, whose exact values are not important and
can  change from one appearance to another. The labeling of the
constants $c_1, c_2, \dots$ starts anew in the statement of each
result. The dependence of the constant $c$ on the dimension $d$ will
not be mentioned explicitly. We will use ``$:=$" to denote a
definition, which is read as ``is defined to be".  We
will use $dx$ to denote the Lebesgue measure in $\bR^d$.
For a Borel set $A\subset \bR^d$, we also use $|A|$ to denote its Lebesgue measure.
We denote the Euclidean distance between
$x$ and $y$ in $\R^d$ by  $|x-y|$ and
denote by $B(x, r)$ the open
ball centered at $x\in \bR^d$ with radius $r>0$.
For $a, b\in \bR$, $a\wedge b:=\min \{a,
b\}$ and $a\vee b:=\max\{a, b\}$.
For any two positive functions $f$ and $g$,
we use the notation $f(r)\asymp
g(r),\ r\to a$ to denote that $\tfrac{f(r)}{g(r)}$ stays between two positive
constants as $r\to a$.
$f\asymp
g$ simply
means that there is a positive constant $c\geq 1$
so that $c^{-1}\, g \leq f \leq c\, g$ on their common domain of
definition.
For any open $D\subset\bR^d$ and $x\in D$,
$\delta_D(x)$ stands for the Euclidean distance between
$x$ and $D^c$.

\section{Scaling conditions and consequences}

Recall that a function $\phi:(0,\infty)\to (0,\infty)$ is a Bernstein function if
it has the  representation
$$
\phi(\lambda)=a+b\lambda + \int_{(0,\infty)}(1-e^{-\lambda t})\, \mu(dt)\, ,
$$
where $a,b\ge 0$ and $\mu$ is a measure on $(0,\infty)$ satisfying $\int_{(0,\infty)}(1\wedge t)\, \mu(dt) <\infty$.
A function $\phi:(0,\infty)\to (0,\infty)$ is a Bernstein function if and only if it is the Laplace exponent of a (killed) subordinator $S=(S_t)_{t\ge 0}$: $\E[\exp \{-\lambda S_t\}]=\exp\{-t \phi(\lambda)\}$ for all $t\ge 0$ and $\lambda>0$.

It is well-known that, if $\phi$ is a Bernstein function, then
\begin{equation}\label{e:Berall}
\phi(\lambda t)\le \lambda\phi(t) \qquad \text{ for all } \lambda \ge 1, t >0\, ,
\end{equation}
implying
\begin{equation}\label{e:uv}
\frac{\phi(v)}{v}\le \frac{\phi(u)}{u}\, ,\quad 0<u\le v\, .
\end{equation}
Note that  \eqref{e:uv} implies
\begin{equation}\label{e:Berall1}
\lambda
\phi'(\lambda)\le
\phi(\lambda)
\qquad \text{ for all }
\lambda >0.
\end{equation}

We remark that, since $\phi$ is increasing, \eqref{e:Berall} is equivalent to that $\phi$
is an $O$-regularly varying function, see \cite[Section 2.0.2]{BGT}.

Clearly \eqref{e:Berall} implies the following observation.

\begin{lemma}\label{l:phi-property}
If $\phi$ is a Bernstein function, then for all  $\lambda, t>0$,
$1 \wedge  \lambda\le {\phi(\lambda t)}/{\phi(t)} \le 1 \vee
 \lambda$.
\end{lemma}

Note that with this lemma, we can replace expressions of the type $\phi(\lambda t)$,
when $\phi$ is a Bernstein function,
with $\lambda >0$ fixed and $t>0$ arbitrary, by $\phi(t)$ up to a multiplicative
constant depending on $\lambda$.
We will often do this without explicitly mentioning it.

In the remainder of this paper, we will always assume that $\phi$ is a complete Bernstein
function, that is, the L\'evy measure $\mu$ of $\phi$ has a completely monotone density.
We will denote this density by $\mu(t)$.
For properties of complete Bernstein function, we refer our reader to \cite{SSV}.

We will assume that $\phi$ satisfies either {\bf (H1)}, or {\bf (H2)}, or both.
Note that it follows from the right-hand side inequality in \eqref{e:H1} that $\phi$ has no drift, i.e., $b=0$.
It also follows from the left-hand side inequality in \eqref{e:H2}
that $\phi$ has no killing term, i.e., $a=0$.
Since for most of this paper we assume both {\bf (H1)} and {\bf (H2)},
it is harmless to immediately assume that $a=b=0$
(regardless whether the scaling conditions hold).
So, from now on, $a=b=0$.

Throughout this paper, we use $S=(S_t)_{t\ge 0}$ to denote a subordinator
with Laplace exponent $\phi$. Since $\phi$ is a complete
Bernstein function, the potential measure $U$ of $S$ has a complete monotone density $u(t)$
(see \cite[Theorem 10.3]{SSV} or \cite[Corollary 13.2.3]{KSV3}),
called the potential density of $S$.

Without loss of generality we further assume that $\phi(1)=1$.
Then by taking $r=1$ and $R=\lambda$ in {\bf (H1)}, and $R=1$ and $r=\lambda$ in {\bf (H2)}, we get that
\begin{equation}\label{e:H1-large}
a_1\lambda^{\delta_1} \le \phi(\lambda ) \le a_2 \lambda^{\delta_2}\, ,
\quad \lambda \ge 1\, ,
\end{equation}
and
\begin{equation}\label{e:H2-small}
a^{-1}_4\lambda^{\delta_4} \le \phi(\lambda ) \le a_3^{-1} \lambda^{\delta_3}\, , \quad \lambda \le 1\, .
\end{equation}
If $0<r<1<R$, using \eqref{e:H1-large} and \eqref{e:H2-small}, we have
that under {\bf (H1)}--{\bf (H2)},
$$
\frac{\phi(R)}{\phi(r)} \le a_2a_4 \frac{R^{\delta_2}}{r^{\delta_4}}  \le a_2a_4  \left(\frac{R}{r}\right)^{\delta_2 \vee \delta_4}
\quad \text{ and }\quad
\frac{\phi(R)}{\phi(r)} \ge a_1a_3 \frac{R^{\delta_1}}{r^{\delta_3}}
\ge a_1a_3  \left(\frac{R}{r}\right)^{\delta_1 \wedge \delta_3}.
$$
Combining these with {\bf (H1)} and {\bf (H2)} we get
\begin{equation}\label{e:sc1}
a_5 \left(\frac{R}{r}\right)^{\delta_1 \wedge \delta_3} \le \frac{\phi(R)}{\phi(r)} \le a_6\left(\frac{R}{r}\right)^{\delta_2 \vee \delta_4}, \quad 0<r<R<\infty\, .
\end{equation}

For $a>0$, we define $\phi^a(\lambda)= \phi(\lambda a^{-2})/\phi( a^{-2})$. Then $\phi^a$
is again a complete Bernstein function satisfying $\phi^a(1)=1$.
We will use $\mu^a(dt)$ and $\mu^a(t)$ to denote the L\'evy measure and L\'evy density
of $\phi^a$ respectively, $S^a=(S^a_t)_{t\ge 0}$ to denote a subordinator with Laplace exponent
$\phi^a$, and $u^a(t)$ to denote the potential density of $S^a$.
Since
$$
\phi^a(\lambda)=\int_0^{\infty}(1-e^{-\lambda t})\, \mu^a(t)dt,
\quad \int^\infty_0e^{-\lambda t}u^a(t)dt=\frac1{\phi^a(\lambda)}, \qquad \lambda>0,
$$
it is straightforward to see that
\begin{equation}\label{e:levy-density-a}
\mu^a(t)=\frac{a^2}{\phi(a^{-2})}\mu(a^2 t)\, , \qquad t>0\, ,
\end{equation}
and
\begin{equation}\label{e:potential-density-a}
u^a(t)=a^2 \phi(a^{-2})u(a^2 t)\, , \qquad t>0\, .
\end{equation}

Now applying this to $\phi^a$, we get
that under {\bf (H1)}--{\bf (H2)},
\begin{equation}\label{e:sc2}
a_5 \left(\frac{R}{r}\right)^{\delta_1 \wedge \delta_3} \le \frac{\phi^a(R)}{\phi^a(r)} \le a_6 \left(\frac{R}{r}\right)^{\delta_2 \vee \delta_4}, \quad a>0,\ 0<r<R<\infty\, .
\end{equation}

The results in the next lemma will be used many times later in the paper.

\begin{lemma}\label{l:integral-estimates-phi}
 Assume {\bf (H1)} and {\bf (H2)}.
    There exists $c=c(a_1, a_2, a_3, a_4, \delta_1, \delta_2, \delta_3, \delta_4)>0$ such that
    \begin{eqnarray}
    \int_0^{\lambda^{-1}} \phi (r^{-2})^{1/2}\, dr    \le  c  \lambda^{-1}\phi(\lambda^{2})^{1/2},  &\quad &\textrm{for all }\lambda > 0\, ,\label{e:ie-1} \\
    \lambda^2 \int_0^{\lambda^{-1}} r \phi (r^{-2})\, dr + \int_{\lambda^{-1}}^\infty r^{-1} \phi (r^{-2})\, dr \le  c \phi(\lambda^2)\, ,&\quad &\textrm{for all } \lambda >0\, ,\label{e:ie-2}
    \\
    \int_0^{\lambda^{-1}} r^{-1}\phi(r^{-2})^{-1}\, dr \le c \phi(\lambda^2)^{-1}\, , &\quad &\textrm{for all } \lambda >0\,.\label{e:ie-3}
    \end{eqnarray}
\end{lemma}
\pf
This result is essentially Karamata's theorem for $O$-regularly varying functions with constants controlled
and its proof is hidden in the proofs in \cite[Section 2.6]{BGT}.
Taking into account \eqref{e:sc1}, direct proofs of
\eqref{e:ie-1}--\eqref{e:ie-3}
are the same as those of \cite[Lemma 4.1]{KSV7}.
We omit the proof here. \qed

The following result plays a crucial role in this paper.

\begin{prop}\label{p:zahle}
Suppose that $w$ is a completely monotone function given by
$$
w(t)=\int^\infty_0e^{-st} f(s)\, ds,
$$
where $f$ is a nonnegative decreasing function.
\begin{itemize}
    \item[(a)] It holds that
    \begin{equation}\label{e:zahle-00}
    f(s)\le \left(1-e^{-1}\right)^{-1} s^{-1}w(s^{-1})\, , \quad s>0.
        \end{equation}
    \item[(b)] If there exist $\delta\in (0, 1)$ and $a, s_0>0$ such that
    \begin{equation}\label{e:zahle-1}
    w(\lambda t)\le a \lambda^{-\delta} w(t), \quad \lambda\ge 1, t\ge 1/s_0,
    \end{equation}
    then there exists $c_1=c_1(w,a,s_0, \delta)>0$ such that
    $$
    f(s)\ge c_1 s^{-1}w(s^{-1}), \quad s\le s_0.
    $$
    \item[(c)] If there exist $\delta\in (0, 1)$ and $a, s_0>0$ such that
    \begin{equation}\label{e:zahle-2}
    w(\lambda t)\ge a \lambda^{-\delta} w(t), \quad \text{for all }  \lambda   \le 1 \text{ and } t\le 1/s_0,
    \end{equation}
    then there exists $c_2=c_2(w,a,s_0, \delta)>0$  such that
    $$
    f(s)\ge c_2 s^{-1}w(s^{-1}), \quad s\ge s_0.
    $$
\end{itemize}
\end{prop}

\pf
This result follows from Karamata's Tauberian theorem and monotone density theorem (together
with their counterparts at 0) for $O$-regularly varying functions, see \cite[Theorem 2.10.2 and
Proposition 2.10.3]{BGT}. Here we give a direct proof.

Direct proofs of (a) and (b) are given in \cite{Z}  (see also \cite[Proposition 13.2.5]{KSV3}).

\noindent (c)
Let $\rho:=\left(\int_0^{s_0}e^{-s}f(s)\,d s\right)\left(\int_{s_0}^{\infty}e^{-s} f(s)\, ds\right)^{-1}$.
Note that $\rho=\rho(f,s_0)=\rho(w, s_0)$.
For any $t\le 1$, we have
\begin{eqnarray*}
&&\int_0^{s_0}e^{-t s}\, f(s)\, ds
=\int_0^{s_0}e^{(1-t) s}e^{-s}\, f(s)\, ds
\le e^{(1-t)s_0}\int_0^{s_0}e^{-s}\, f(s)\, ds\\
& &=  \rho e^{(1-t)s_0}\int_{s_0}^\infty e^{-s}\, f(s)\, ds\,\le\,  \rho \int_{s_0}^\infty e^{-t s}\, f(s)\, ds.
\end{eqnarray*}
Thus for any $t \le 1$
\begin{eqnarray*}
w(t)\le (\rho +1) \int^\infty_{s_0} e^{-st} f(s)\, ds
=\frac{\rho+1}{t}\int^\infty_{s_0 t} e^{-s}f\left(\frac{s}{t}\right)\, ds.
\end{eqnarray*}

Let $t\le 1$ be arbitrary. For any $r\in (0, 1]$,  we have
\begin{eqnarray}
tw(t)& \le & (\rho+1)\int^r_0 {\bf 1}_{\{ s_0 t <s\}}e^{-s}f\left(\frac{s}{t}\right)\, ds + (\rho+1)\int_r^\infty
e^{-s}{\bf 1}_{\{ s_0 t <s\}} f\left(\frac{s}{t}\right)ds\nn\\
&\le& (\rho+1) \int^r_0 {\bf 1}_{\{s_0 t <s\}} e^{-s}f\left(\frac{s}{t}\right)\, ds
+(\rho +1) f\left(\frac{r}{t}\right)e^{-r}\nn\\
&\le& (\rho+1)\left(1-e^{-1}\right)^{-1}t\int^r_0{\bf 1}_{\{ s_0 t <s\}}e^{-s}\frac1s\, w\left(\frac{t}{s}\right)\, ds
+(\rho +1) f\left(\frac{r}{t}\right)e^{-r}, \label{e:n123}
\end{eqnarray}
where in the last line we used
\eqref{e:zahle-00}.

Now we assume \eqref{e:zahle-2} and apply it to $w(\tfrac{t}{s})$ in \eqref{e:n123}.
Note that $s\le r\le 1$, and since $s_0 t <s$ we also have that $t\le s/s_0$. Thus
$
w(\tfrac{t}{s})\le a^{-1} s^{\delta} w(t)\, ,
$
implying that
$$
t w(t)\le (\rho+1)a^{-1} \left(1-e^{-1}\right)^{-1} t w(t) \int^r_0{\bf 1}_{\{ s_0 t <s\}}e^{-s}s^{\delta-1}\, ds
+ (\rho +1) f\left(\frac{r}{t}\right)e^{-r}\, .
$$
Choose $r=r(a,s_0,\delta)\in(0, 1]$ small enough so that
$$
(\rho +1) a^{-1} \left(1-e^{-1}\right)^{-1}\int_0^r e^{-s}s^{\delta-1}\, ds \le \frac12.
$$
For this choice of $r$, we have
$
f(\tfrac{r}{t})\ge c_1 tw(t)$,  $t\le 1$,
for some $c_1=c_1(a,w,a,s_0 )>0$. Thus
$$
f(s)\ge c_1\frac{r}{s}w\left( \frac{r}{s}\right) \ge c_2 s^{-1}w(s^{-1}), \quad s\ge r,
$$
where $c_2=c_1r$. In order to extend the inequality to $s\ge s_0$ it suffices to use the continuity of $w$.
\qed

\begin{corollary}\label{c:asymptotics-u}
\begin{itemize}
    \item[(a)] The potential density $u$ of $S$ satisfies
    \begin{equation}\label{e:u-upper-bound}
    u(t)\le (1-e^{-1})^{-1}  t^{-1}\phi(t^{-1})^{-1}\, ,\quad t>0.
    \end{equation}
    \item[(b)] If {\bf (H1)} holds, then there exists $c_1=c_1(\phi)>0$ such that
    \begin{equation}\label{e:u-lower-bound-zero}
    u(t)\ge c_1  t^{-1}\phi(t^{-1})^{-1}\, ,\quad 0<t\le 1.
    \end{equation}
    \item[(c)] If {\bf (H2)} holds, then there exists $c_2=c_2(\phi)>0$ such that
    \begin{equation}\label{e:u-lower-bound-infty}
    u(t)\ge c_2  t^{-1}\phi(t^{-1})^{-1}\, ,\quad 1\le t <\infty.
    \end{equation}
\end{itemize}
\end{corollary}
\pf
(a) The claim follows from Proposition \ref{p:zahle}(a) with
$
w(t):=\int_0^{\infty}e^{-s t}u(s)\, ds =\tfrac{1}{\phi(t)}\, .
$

\noindent
(b) By the left-hand side of \eqref{e:H1}, $w(t)=\phi(t)^{-1}$ satisfies \eqref{e:zahle-1} with $\delta=\delta_1$, $a=a_1^{-1}$ and $s_0=1$. The claim follows from Proposition
\ref{p:zahle}(b) with $c_1=c_1(\phi, a_1, \delta_1)$.

\noindent
(c)
By the left-hand side of \eqref{e:H2}, $w(t)=\phi(t)^{-1}$
satisfies \eqref{e:zahle-2} with $\delta=\delta_3$, $a=a_3^{-1}$ and $s_0=1$.
The claim follows from Proposition \ref{p:zahle}(c)
with $c_2=c_2(\phi, a_3, \delta_3)$.
\qed

Since $\phi$ is a complete Bernstein function, its conjugate function $\phi^{\ast}(\lambda):=\frac{\lambda}{\phi(\lambda)}$ is also complete Bernstein. It is immediate to see that,
under {\bf (H2)} for $\phi$, the function $\phi^{\ast}$ satisfies
\begin{align*}
a_4^{-1}\left(\frac{R}{r}\right)^{1-\delta_4}\le \frac{\phi(R)}{\phi(r)}\le a_3^{-1}
\left(\frac{R}{r}\right)^{1-\delta_3}, \quad r\le R\le 1.
\end{align*}
Since the potential density $u^{\ast}$ of $\phi^{\ast}$ is equal to the tail $\mu(t,\infty)$ of the L\'evy measure $\mu$ (see \cite[Corollary 5.5]{BBKRSV}), we conclude from Corollary \ref{c:asymptotics-u} that
\begin{align}
& \mu(t,\infty) \le (1-e^{-1})^{-1} t^{-1} \phi^{\ast}(t^{-1})^{-1}\, ,\quad t>0\, ,\label{e:tail-upper-bound}\\
& \mu(t,\infty) \ge c t^{-1} \phi^{\ast}(t^{-1})^{-1}\, ,\quad 1\le t <\infty\, ,\quad \textrm{if {\bf (H2)} holds}\, .
\label{e:tail-lower-bound-infty}
\end{align}

\begin{prop}\label{p:asymptotics-mu}
\begin{itemize}
        \item[(a)] The L\'evy density $\mu$ of $S$ satisfies
    \begin{equation}\label{e:mu-upper-bound}
    \mu(t)\le (1-2e^{-1})^{-1}  t^{-1}\phi(t^{-1})\, ,\quad t>0.
    \end{equation}
    \item[(b)] If {\bf (H1)} holds, then there exists $c_1=c_1(\phi)>0$ such that
    \begin{equation}\label{e:mu-lower-bound-zero}
    \mu(t)\ge c_1  t^{-1}\phi(t^{-1})\, ,\quad 0<t\le 1.
    \end{equation}
    \item[(c)] If {\bf (H2)} holds,  then there exists $c_2=c_2(\phi)>0$ such that
    \begin{equation}\label{e:mu-lower-bound-infty}
    \mu(t)\ge c_2  t^{-1}\phi(t^{-1})\, ,\quad 1\le t <\infty.
    \end{equation}
\end{itemize}
\end{prop}
\pf (a) This is proved in \cite[Lemma A.1, Proposition 3.3]{KM}.

\noindent
(b) This is proved in \cite[Theorem 13.2.10]{KSV3}.

\noindent
(c) The proof is similar to the proof of (b).
It follows from \eqref{e:tail-upper-bound} and \eqref{e:tail-lower-bound-infty} that there exists a constant $c_1>0$ such that
$
c_1^{-1}\phi(s^{-1})\le u^{\ast}(s) \le c_1 \phi(s^{-1})$ for $1\le s <\infty$.
Fix $\lambda:=(2c_1^2 a_3^{-1})^{1/\delta_3}\vee1\ge1$.
Then by the left-hand side of {\bf (H2)}, we have that for $s\ge \lambda$,
$$
u^{\ast}(s)\le c_1\phi(s^{-1})=c_1\phi(\lambda^{-1}(\lambda^{-1} s)^{-1})
\le c_1 a_3^{-1} \lambda^{-\delta_3}\phi((\lambda^{-1}s)^{-1})
\le c_1^2 a_3^{-1}\lambda^{-\delta_3} u^{\ast}(\lambda^{-1}s) \le\frac12 u^{\ast}(\lambda^{-1}s)
$$
by our choice of $\lambda$. Further,
$$
(1-\lambda^{-1})s\mu(\lambda^{-1}s)\ge \int^s_{\lambda^{-1} s}\mu(t)\, dt=u^{\ast}(\lambda^{-1} s)-u^{\ast}(s)\ge u^{\ast}(\lambda^{-1} s) -\frac12 u^{\ast}(\lambda^{-1} s)=\frac12 u^{\ast}(\lambda^{-1} s)\, .
$$
This implies that for all
$t\ge 1$
$$
\mu(t)\ge \frac{1}{2(1-\lambda^{-1})\lambda} t^{-1} u^{\ast}(t)= c_2 t^{-1}u^{\ast}(t)\ge c_3 t^{-1}\phi(t^{-1})
$$
for some constants $c_2, c_3>0$.
\qed

We conclude this section with some conditions on $\phi$ which imply {\bf (H1)} and {\bf (H2)}.

\noindent
{\bf ($H_0$):} There exist $\beta\in (0,2)$ and a function $\ell:(0,\infty)\to (0,\infty)$ which is measurable,  bounded on compact subsets of $(0,\infty)$  and slowly varying at 0 such that
\begin{equation}\label{e:reg-var}
\phi(\lambda) \asymp \lambda^{\beta/2}\ell(\lambda)\, ,\quad \lambda \to 0+\, .
\end{equation}

\noindent {\bf ($H_{\infty}$):} There exist $\alpha\in (0,2)$ and a function $\widetilde{\ell}:(0,\infty)\to (0,\infty)$ which is measurable,  bounded on compact subsets of $(0,\infty)$ and slowly varying at infinity such that
\begin{equation}\label{e:reg-var-infty}
\phi(\lambda) \asymp \lambda^{\alpha/2}\widetilde{\ell}(\lambda)\, ,\quad \lambda \to \infty\, .
\end{equation}
Using Potter's theorem (cf.~\cite[Theorem 1.5.6]{BGT}), it is proved in \cite{KSV3} that {\bf ($H_{\infty}$)} implies the right-hand side inequality of {\bf (H1)}.
One can similarly prove that {\bf ($H_{\infty}$)} also implies the left-hand side inequality of {\bf (H1)}
and that {\bf ($H_0$)} implies {\bf (H2)}.


\section{Applications to subordinate Brownian motions}

Recall that $S=(S_t)_{t\ge 0}$ is a subordinator with Laplace exponent $\phi$.
Let $W=(W_t, \P_x)_{t\ge 0}$ be a $d$-dimensional Brownian motion independent of $S$ and
with transition density
$$
q(t,x,y)=(4\pi t)^{-d/2} e^{-\frac{|x-y|^2}{4t}}\, ,\quad x,y\in \R^d, \ t>0\, .
$$
The process $X=(X_t)_{t\ge 0}$ defined by $X_t:=W(S_t)$ is called a subordinate Brownian motion.
It is a rotationally invariant L\'evy process with characteristic exponent
$\phi(|\xi|^2)$, $\xi\in \R^d$,
and transition density given by
$$
p(t, x, y)=\int^\infty_0q(s, x, y)\P(S_t\in ds).
$$
By spatial homogeneity, the L\'evy measure of $X$
 has a density $J(x)=j(|x|)$, where $j:(0,\infty)\to (0,\infty)$ is given by
\begin{equation}\label{e:jumping-function}
j(r):=\int_0^{\infty}(4\pi t)^{-d/2} e^{-r^2/(4t)}\mu(t)\, dt\, .
\end{equation}
Note that $j$ is continuous and decreasing.
We define $J(x, y):=J(y-x)$.

By the Chung-Fuchs criterion the process $X$ is transient if and only if
\begin{equation}\label{e:CF}
\int_0^1 \frac{\lambda^{d/2-1}}{\phi(\lambda)}\, d\lambda <\infty\, .
\end{equation}
Note that if $d\ge 3$, then  $X$ is always transient.
If {\bf (H2)} holds and $d>2 \delta_4$, then $X$ is transient. In particular,
if {\bf (H2)} holds and $d\ge 2$, then $X$ is transient.
When $X$ is transient, the mean occupation time measure of $X$ admits a density $G(x,y)=g(|x-y|)$
which is called the Green function of $X$,
and is given by the formula
\begin{equation}\label{e:green-function}
g(r):=\int_0^{\infty}(4\pi t)^{-d/2} e^{-r^2/(4t)}u(t)\, dt\, .
\end{equation}
Here $u$ is the potential density of the subordinator $S$.
Note that by the transience assumption, the integral converges. Moreover, $g$ is continuous and decreasing.

We first record the upper bounds of $j(r)$ and $g(r)$.
\begin{lemma}\label{l:j-g-upper}
 \item[(a)] It holds that
$
    j(r)\le c_1 r^{-d}\phi(r^{-2})$ for all $r > 0.
 $
    \item[(b)] If $d\ge 3$  then
   $
    g(r)\le c_2 r^{-d}\phi(r^{-2})^{-1}$ for all $r > 0.
   $
\end{lemma}
\pf
(a) We write
$$
j(r)=\int_0^{r^2}(4\pi t)^{-d/2} e^{-r^2/(4t)}\mu(t)\, dt + \int_{r^2}^{\infty}(4\pi t)^{-d/2} e^{-r^2/(4t)}\mu(t)\, dt:=
J_1 + J_2\, .
$$
To estimate $J_2$ from above we first use \eqref{e:mu-upper-bound} and then the monotonicity of $\phi$ to obtain
\begin{eqnarray*}
J_2&\le & c_1 \int_{r^2}^{\infty}(4\pi t)^{-d/2} e^{-r^2/(4t)} t^{-1}\phi(t^{-1})\, dt \\
&\le &c_1 \phi(r^{-2}) \int_0^{\infty}t^{-d/2-1}e^{-r^2/(4t)} \, dt \,=\,  c_2 \phi(r^{-2}) r^{-d}\, .
\end{eqnarray*}
By \eqref{e:uv} we have $\phi(t^{-1})/t^{-1} \le \phi(r^{-2})/r^{-2}$ for $t\le r^2$ (i.e., $r^{-2}\le t^{-1}$),
thus
by \eqref{e:mu-upper-bound}
\begin{eqnarray*}
&&J_1\le  c_3 \int_0^{r^2} t^{-d/2} e^{-r^2/(4t)} t^{-1}\phi(t^{-1})\, dt \,\le \,
c_3 r^2\phi(r^{-2}) \int_0^{r^2} t^{-d/2-2} e^{-r^2/(4t)}\, dt  \\
&&\le  c_3 r^2\phi(r^{-2}) \int_0^{\infty} t^{-d/2-2} e^{-r^2/(4t)}\, dt  \, = \, c_4 r^{-d}\phi(r^{-2}).
\end{eqnarray*}

\noindent
(b) We write
\begin{eqnarray}
g(r)=\int_0^{r^2}(4\pi t)^{-d/2} e^{-r^2/(4t)}u(t)\, dt + \int_{r^2}^{\infty}(4\pi t)^{-d/2}
e^{-r^2/(4t)}u(t)\, dt:=L_1+L_2\, . \label{e:gert}
\end{eqnarray}
By using \eqref{e:u-upper-bound} in the first inequality and
the monotonicity of $\phi$ in the second inequality, we get
\begin{eqnarray}
L_1 &\le  & c_5 \int_0^{r^2} (4\pi t)^{-d/2} e^{-r^2/(4t)} t^{-1} \phi(t^{-1})^{-1}\, dt
\le c_6 \phi(r^{-2})^{-1} \int_0^{r^2} t^{-d/2-1} e^{-r^2/(4t)}\, dt\nonumber\\
& \le & c_6 \phi(r^{-2})^{-1} \int_0^{\infty} t^{-d/2-1} e^{-r^2/(4t)}\, dt
\,= \, c_7\phi(r^{-2})^{-1}r^{-d}\, . \label{e:df3}
\end{eqnarray}
Since $d\ge 3$, using that $u$ is decreasing in the second inequality and
\eqref{e:u-upper-bound} in the third, we get
\begin{eqnarray*}
L_2 \,\le\,   c_8 \int_{r^2}^{\infty}t^{-d/2}u(t)\, dt\,\le\,  c_8 u(r^2)
\int_{r^2}^{\infty}t^{-d/2}\, dt\,
\le\,  c_9 \phi(r^{-2})^{-1} r^{-d}\, .
\end{eqnarray*}
\qed

Our next goal is to establish the asymptotic behaviors of $j(r)$ and $g(r)$ for small and/or large $r$
under {\bf (H1)} or {\bf (H2)}, or both.

\begin{lemma}\label{l:j-g-near-origin} Assume {\bf (H1)}.
\begin{itemize}
    \item[(a)] It holds that
    \begin{equation}\label{e:j-near-origin}
    j(r)\asymp r^{-d}\phi(r^{-2})\, ,\qquad r\to 0.
    \end{equation}
    \item[(b)] If $d>2\delta_2$ and $X$ is transient, then
    \begin{equation}\label{e:g-near-origin}
    g(r)\asymp r^{-d}\phi(r^{-2})^{-1}\, ,\qquad r\to 0.
    \end{equation}
\end{itemize}
\end{lemma}
\pf
(a) is proved in \cite[Theorem 13.3.2]{KSV3}, so we only prove (b).
First note that the assumption $d> 2\delta_2$ is always satisfied when $d\ge 2$.

By Lemma \ref{l:j-g-upper}
we only need to prove the upper bound in \eqref{e:g-near-origin} for $d \le 2$.
To do that  we write
$g(r)=L_1+L_2$
as in \eqref{e:gert}.
First note that, by the same argument as for \eqref{e:df3}, we have
$L_1 \le   c_1\phi(r^{-2})^{-1}r^{-d}.$

Let $d\le 2$ and $r \le 1$. We split $L_2$ into two parts:
$$
L_2 \le c_2\int_{r^2}^1 t^{-d/2}  u(t)\, dt+ c_2\int_1^{\infty} t^{-d/2} e^{-r^2/(4t)} u(t)\, dt
=: L_{21}+L_{22}\, .
$$
For $L_{21}$ we use \eqref{e:u-upper-bound} and the change of variables $t=r^2 s$ to get
\begin{eqnarray*}
L_{21}
= c_3 r^{-d}\int_1^{r^{-2}}s^{-d/2-1}\phi(r^{-2}s^{-1})^{-1}\, ds\, .
\end{eqnarray*}
Since $0<r\le 1$ and $r^2\le s^{-1}\le 1$, it follows from
\eqref{e:H1}
 that $\phi(r^{-2}s^{-1})^{-1}\le a_2 s^{\delta_2}\phi(r^{-2})^{-1}$. Hence
$$
L_{21}\le c_4 r^{-d}\phi(r^{-2})^{-1}\int_1^{\infty}s^{-d/2-1+\delta_2}\, ds =c_5 r^{-d}\phi(r^{-2})^{-1}\, ,
$$
since the integral converges under the assumption $d>2\delta_2$.
Note that using  ({\bf H1}) and  the assumption that $2\delta_2 < d$, we have
$r^{-d}\phi(r^{-2})^{-1} \ge c_{6}r^{2\delta_2-d}\ge c_{6}>0$.
Since  $L_{22}$ is bounded for $r\le 1$ by \cite[Lemma 4.4]{KM}, we have proved the upper bound.

To prove the converse inequality for all $d \ge 1$, we use \eqref{e:u-lower-bound-zero} in the second inequality and \eqref{e:uv} in the third to get that  for $r \le 1$,
\begin{align*}
&g(r)\ge  (4\pi)^{-d/2}\int_0^{1/r^2}(t r^2)^{-d/2}e^{-r^2/(4tr^2)} u(r^2 t) r^2\, dt\\
&\ge   c_7 r^{2-d}\int_0^1 t^{-d/2}e^{-1/(4t)} r^{-2}t^{-1} \phi(r^{-2}t^{-1})^{-1}\, dt \,\ge\,  c_{8} r^{-d}\phi(r^{-2})^{-1}\, .
\end{align*}
\qed

\begin{lemma}\label{l:j-g-near-infty} Assume {\bf (H2)}.
\begin{itemize}
    \item[(a)] It holds that
    \begin{equation}\label{e:j-near-infty}
    j(r)\asymp r^{-d}\phi(r^{-2})\, ,\qquad r\to \infty.
    \end{equation}
    \item[(b)] If $d>2\delta_4$, then $X$ is transient and
    \begin{equation}\label{e:g-near-infty}
    g(r)\asymp r^{-d}\phi(r^{-2})^{-1}\, ,\qquad r\to \infty.
    \end{equation}
\end{itemize}
\end{lemma}
\pf (a)
By Lemma \ref{l:j-g-upper}
we only need to prove the lower bound in \eqref{e:j-near-infty}.
For the lower bound we have
\begin{eqnarray*}
j(r)&\ge &
(4\pi)^{-d/2} \int_0^1 (r^2 t)^{-d/2} e^{-1/(4t)} \mu(r^2 t) r^2 dt \\
& \ge & c_{1} r^{-d+2}\mu(r^2) \int_0^1 t^{-d/2}e^{-1/(4t)}\, dt
\,\ge
\,
c_{2} r^{-d} \phi(r^{-2})\, ,
\end{eqnarray*}
where in the last inequality we used \eqref{e:mu-lower-bound-infty}.

\noindent
(b)
By \eqref{e:H2-small}, $a_4^{-1}\lambda^{\delta_4}\le \phi(\lambda)$ for all $\lambda \le 1$,
so using the assumption $d>2\delta_4$,
$X$ is transient by \eqref{e:CF}.
Let $r\ge 1$. By the change of variables $s=r^2/t$ we get that
\begin{eqnarray*}
g(r)
=c_3 r^{-d+2}\int_0^{\infty}s^{d/2-2}e^{-s/4} u(r^2 s^{-1})\, ds.
\end{eqnarray*}
By \eqref{e:u-upper-bound}, we have $u(r^2 s^{-1})\le c_4 r^{-2}s \phi(r^{-2}s)^{-1}$. Hence
\begin{eqnarray*}
g(r)\le c_5 r^{-d}\int_0^1s^{d/2-1} \phi(r^{-2}s)^{-1}\, ds+c_5 r^{-d}\int_1^{\infty}s^{d/2-1}e^{-s/4}
\phi(r^{-2}s)^{-1}\, ds\,=:\, L_1+L_2\, .
\end{eqnarray*}
To estimate $L_1$ from above, we note that, by ({\bf H2}),
we have $\phi(r^{-2}s)\ge a_4^{-1} s^{\delta_4}\phi(r^{-2})$, $0<s\le 1$. Hence
$$
L_1 \le  c_6 r^{-d} \phi(r^{-2})^{-1} \int_0^1 s^{d/2-1-\delta_4}e^{-s/4}\, ds
=  c_7 r^{-d}\phi(r^{-2})^{-1}
$$
since the integral converges under the assumption $d/2>\delta_4$.
In order to estimate $L_2$, we use that $r^{-2}\le r^{-2}s$ for $s\ge 1$,
hence by the monotonicity of $\phi$, $\phi(r^{-2}s)\ge \phi(r^{-2})$. Therefore,
$$
L_2 \le c_8 r^{-d} \phi(r^{-2})^{-1}\int_1^{\infty}s^{d/2-1}e^{-s/4}\, ds =c_9 r^{-d} \phi(r^{-2})^{-1}\, .
$$

For the lower bound we have
\begin{eqnarray*}
g(r)&\ge &
(4\pi)^{-d/2} \int_0^1 (r^2 t)^{-d/2} e^{-1/(4t)} u(r^2 t) r^2 dt \\
& \ge & c_{10} r^{-d+2}u(r^2) \int_0^1 t^{-d/2}e^{-1/(4t)}\, dt
\,\ge
\,
c_{11} r^{-d} \phi(r^{-2})^{-1}\, ,
\end{eqnarray*}
where in the last inequality we used the left inequality in \eqref{e:u-lower-bound-infty}.
\qed

We now have the asymptotic behaviors of the Green function and L\'evy density of $X$ as
 an immediate consequence of Lemmas \ref{l:j-g-near-origin}--\ref{l:j-g-near-infty}.

\begin{thm}\label{t:J-G}
Assume  both {\bf (H1)} and {\bf (H2)}.
\begin{itemize}
    \item[(a)] It holds that
    \begin{equation}\label{e:J}
    J(x)\asymp |x|^{-d}\phi(|x|^{-2})\, , \quad \textrm{for all }x\neq 0\, .
    \end{equation}
    \item[(b)] If $d > 2(\delta_2\vee \delta_4)$ then the process $X$ is transient  and it holds
    \begin{equation}\label{e:G}
    G(x)\asymp |x|^{-d}\phi(|x|^{-2})^{-1}\, , \quad \textrm{for all }x\neq 0\, .
    \end{equation}
\end{itemize}
\end{thm}

We record a simple consequence of Theorem \ref{t:J-G}.
\begin{corollary}\label{c:doubling-condition}
Assume {\bf (H1)} and {\bf (H2)}.
There exists $c>0$ such that  $J(x)\,\le\, c \,J(2x)$ for all $x\neq 0$ and, if $d > 2(\delta_2\vee \delta_4)$ then $G(x)\,\le\, c\, G(2x)$ for all $x\neq 0$.
\end{corollary}
\pf By Theorem \ref{t:J-G} there exists $c_1>0$ such that
$$
\frac{J(x)}{J(2x)}\le c_1 \frac{|x|^{-d}\phi(|x|^{-2})}{|2x|^{-d}\phi(|2x|^{-2})}=2^d c_1 \frac{\phi(|x|^{-2})}{\phi(4^{-1}|x|^{-2})} \le c_2\, ,\quad x\neq 0\, ,
$$
where the last inequality follows from Lemma \ref{l:phi-property}.
The statement about $G$ is proved in the same way.
\qed

We also record the following property of $j$: There exists $c>0$ such that
\begin{equation}\label{e:j-decay}
j(r)\le cj(r+1)\, ,\quad \textrm{for all }r\ge 1\, .
\end{equation}
This is a consequence of the similar property of $\mu(t)$ and is proved in \cite[Proposition 13.3.5]{KSV3}.
By Corollary \ref{c:doubling-condition} we also have
\begin{equation}\label{e:doubling-condition}
j(r)\le c j(2r)\, ,\quad r>0\, .
\end{equation}

\medskip
Let $a>0$. Recall that $\phi^a$ was defined by $\phi^a(\lambda)=\phi(\lambda a^{-2})/\phi(a^{-2})$. Let $S^a=(S^a)_{t\ge 0}$ be a subordinator with Laplace exponent $\phi^a$ independent of the Brownian motion $W$. Let $X^a=(X^a_t)_{t\ge 0}$ be defined by $X^a_t:=W_{S^a_t}$. Then $X^a$ is a rotationally invariant L\'evy process with characteristic exponent
$$
\phi^a(|\xi|^2)=
\frac{\phi(a^{-2}|\xi|^2)}{\phi(a^{-2})}\, , \quad \xi\in \R^d\, .
$$
This shows that $X^a$ is identical in law to the process $\{a^{-1}X_{t/\phi(a^{-2})}\}_{t\ge 0}$.

The L\'evy measure of $X^a$ has a density $J^a(x)=j^a(|x|)$, where $j^a$ is given by
\begin{eqnarray}
&&j^a(r)=\int_0^{\infty}(4\pi t)^{-d/2} e^{-r^2/(4t)}\mu^a(t)\, dt
=\int_0^{\infty}(4\pi t)^{-d/2} e^{-r^2/(4t)}\frac{a^2}{\phi(a^{-2})}\mu(a^2t)\, dt \nonumber\\
&&=a^d \phi(a^{-2})^{-1} \int_0^{\infty}(4\pi s)^{-d/2}e^{-a^2 r^2/(4s)}\mu(s)\, ds= a^d \phi(a^{-2})^{-1} j(ar)\, . \label{e:jumping-function-a}
\end{eqnarray}
In the second line we used \eqref{e:levy-density-a} and in the third the change of variables $s=a^2t$.
Together with Theorem \ref{t:J-G}(a), \eqref{e:jumping-function-a} gives the following corollary.
\begin{corollary}\label{c:J-G-a}
Assume {\bf (H1)} and {\bf (H2)}.
There exists $c>1$ such that for all $a>0$ and all $x\neq 0$,
\begin{align}
c^{-1}  |x|^{-d}\phi^a(|x|^{-2}) & \le J^a(x) \le c |x|^{-d}\phi^a(|x|^{-2})\, .\label{e:J-a}
 \end{align}
\end{corollary}

Define
$$
\Phi(r):=\frac{1}{\phi(r^{-2})}\, ,\quad r>0\, .
$$
Then $\Phi$ is a strictly increasing function satisfying $\Phi(1)=1$.  In terms of $\Phi$, we can rewrite \eqref{e:J-a} as
\begin{equation}\label{e:J-Phi}
c^{-1}\, \frac{1}{|x|^d\Phi(x)}\le J(x) \le c \frac{1}{|x|^d\Phi(x)}\, .
\end{equation}
Further, \eqref{e:sc2} reads as
\begin{equation}
a_5 \left(\frac{R}{r}\right)^{2(\delta_1 \wedge \delta_3)} \le \frac{\Phi(R)}{\Phi(r)} \le a_6 \left(\frac{R}{r}\right)^{2(\delta_2 \vee \delta_4)}, \quad 0<r<R<\infty\, .\label{e:CK1}
\end{equation}
This implies that
\begin{equation}
\int_0^r\frac {s}{\Phi (s)}\, ds \, \leq \, \frac{a_6}{2(1-\delta_2 \vee \delta_4)} \, \frac{ r^2}{\Phi (r)}\, , \qquad   \textrm{for all } r>0.\label{e:CK2}
\end{equation}
The last three displays show that the process $X$ satisfies conditions (1.4), (1.13) and (1.14) from \cite{CK}. Therefore, by \cite[Theorem 4.12]{CK}, $X$ satisfies the parabolic Harnack inequality, hence also the Harnack inequality.
Thus the following global Harnack inequality is true.
We recall that a function $u:\R^d \to [0,\infty)$ is \emph{harmonic}  with respect to the process $X$ in an open set $D$ if for every relatively compact open set $B\subset D$ it holds that
$$
u(x)=\E_x[u(X_{\tau_B})] \qquad \textrm{ for all } x\in B\, ,
$$
where $\tau_B=\inf\{t>0:\, X_t \notin B\}$ is the exit time of $X$ from $B$.
\begin{thm}\label{t:uhi}
Assume {\bf (H1)} and {\bf (H2)}.
There exists $c=c(\phi)>0$ such that, for any
$r>0$,
$x_0\in \R^d$, and any function $u$ which is nonnegative on $\R^d$ and harmonic with respect to $
X$ in $B(x_0, r)$, we have
$$
u(x)\le c u(y), \quad \textrm{for all }x, y\in B(x_0, r/2).
$$
\end{thm}
This theorem can be also deduced by using the approach in \cite{KSV3}.

We now give some other consequences of \eqref{e:sc2} and Corollary \ref{c:J-G-a}.

Let $B=(B_t, \P_x)_{t\ge 0}$ be a one-dimensional Brownian motion independent of $S^a$
and let $Z^a=(Z^a_t)_{t\ge 0}$ be the one-dimensional subordinate Brownian motion defined by $Z_t:=B(S^a_t)$.
Let $\chi^a$ be the Laplace exponent of the ladder height process of $Z^a$,
$v^a$ be the potential density  of the ladder height process of $Z^a$, and $V^a(t)=\int^t_0v^a(s)ds$ the corresponding renewal function.
It follows from \cite[Corollary
9.7]{Fris} that
$$
\chi^a(\lambda)=
\exp\left(\frac1\pi\int^{\infty}_0\frac{\log(
 \phi^a(\lambda^2\theta^2))}{1+\theta^2}d\theta
\right), \quad \textrm{ for all }
 a,  \lambda>0\, .
$$
Using this and the fact that $\phi^a(\lambda)=\phi(\lambda a^{-2})/\phi(a^{-2})$ we see that
$\chi^a(\lambda)=\phi(a^{-2})^{-1/2} \chi(\lambda/a)$.
This and the identity
$
\int^\infty_0 e^{-\lambda t}v^a(t)\, dt=\tfrac1{\chi^a(\lambda)}
$
imply that  for all $a>0$ and $r>0$, $v^a(t)= a \sqrt{\phi(a^{-2})}v(at)$ so that
\begin{equation}\label{e:Va}
V^a(r)= \int^t_0a \sqrt{\phi(a^{-2})}v(at)ds=  \sqrt{\phi(a^{-2})} V(ar), \quad \text{for all } a, r>0.
\end{equation}
Furthermore, by combining \cite[Proposition 2.6]{KSV5} and \cite[Proposition III.1]{Ber}, we get
\begin{equation}\label{e:V-phi-a}
V^a(r)  \asymp \frac1{\sqrt{\phi^a(r^{-2})}} = \frac{\sqrt{\phi(a^{-2})}}{{\sqrt{\phi(r^{-2} a^{-2})}}}, \quad \text{for all } a, r>0\, .
\end{equation}

\begin{lemma}\label{l:exit-time}
Assume {\bf (H1)} and {\bf (H2)}.
\begin{itemize}
    \item[(a)] There exists $c_1=c_1(\phi)>0$ such that for any  $r>0$ and $x_0 \in \R^d$,
    \begin{eqnarray*}
    \E_x[\tau_{B(x_0,r)}]\le  c_1\, (\phi(r^{-2})\phi((r-|x-x_0|)^{-2}))^{-1/2}\le c\phi(r^{-2})^{-1}\, , \qquad x\in B(x_0, r).
    \end{eqnarray*}
    \item[(b)] There exists $c_2=c_2(\phi)>0$ such that for every $r>0$ and every $x_0\in \R^d$,
    $$
    \inf_{z\in B(x_0,r/2)} \E_z \left[\tau_{B(x_0,r)} \right] \geq c_2 \phi(r^{-2})^{-1}\, .
    $$
\end{itemize}
\end{lemma}

\pf
(a) Using our \eqref{e:V-phi-a} instead of \cite[Proposition 3.2]{KSV7},
the proof of (a) is exactly the same as that of \cite[Lemma 4.4]{KSV7}.

(b) Using \eqref{e:ie-2}, we can repeat the proofs of \cite[Lemmas 13.4.1--13.4.2]{KSV3} to see that
the conclusions of \cite[Lemmas 13.4.1--13.4.2]{KSV3} are valid for all $r>0$. The conclusion of \cite[Lemma 13.4.2]{KSV3}
for all $r>0$ is the desired conclusion in (b).
\qed

The function $J(x,y)$ is the L\'evy intensity of $X$. It determines a L\'evy system for $X$, which describes the jumps of the process
$X$: For any non-negative measurable function $f$ on $\bR_+ \times \bR^d\times \bR^d$ with $f(s, y, y)=0$ for all $y\in \bR^d$, any
stopping time $T$ (with respect to the filtration of $X$) and any $x\in \bR^d$,
\begin{equation}\label{e:levy}
    \E_x \left[\sum_{s\le T} f(s,X_{s-}, X_s) \right]= \E_x \left[ \int_0^T \left( \int_{\bR^d} f(s,X_s, y) J(X_s,y) dy \right) ds \right].
\end{equation}

For every open subset $D\subset \bR^d$, we denote by  $X^D$  the
subprocess of $X$ killed upon exiting $D$.
A subset $D$ of $\R^d$ is said to be Greenian (for $X$) if $X^{D}$ is transient.
For an open Greenian set $D\subset \R^d$, let $G_D(x,y)$ denote the Green function of the killed process $X^D$, and let $K_D(x,z)$ be the Poisson kernel of $D$ with respect to $X$.
Then, by \eqref{e:levy},
\begin{equation}\label{e:poisson-kernel-def}
K_D(x,z)=\int_{\overline{D}^c} G_D(x,y)J(y,z)\, dy\, .
\end{equation}

\begin{prop}\label{p:poisson-kernel-estimate}
Assume {\bf (H1)} and {\bf (H2)}.
There exist $c_1=c_1 (\phi)>0$ and $c_2=c_2 (\phi)>0$ such that for every $r >0$ and $x_0 \in \R^d$,
\begin{eqnarray}
K_{B(x_0,r)}(x,y) \,&\le &\, c_1 \, j(|y-x_0|-r) \left(\phi(r^{-2})\phi((r-|x-x_0|)^{-2})\right)^{-1/2} \label{e:pke-upper}\\
 &\le &\, c_1 \, j(|y-x_0|-r) \phi(r^{-2})^{-1}\ \label{e:pke-upper-worse}
\end{eqnarray}
for all $(x,y) \in B(x_0,r)\times \overline{B(x_0,r)}^c$ and
\begin{equation}\label{e:pke-lower}
K_{B(x_0, r)}(x_0, y) \,\ge\, c_2\, j(|y-x_0|) \phi(r^{-2})^{-1}, \qquad \textrm{ for all } y \in \overline{B(x_0, r)}^c.
\end{equation}
\end{prop}

\pf
With Lemma \ref{l:exit-time} in hand, the proof of this proposition is exactly the same as that of \cite[Lemma 13.4.10]{KSV3}.
\qed

Let $C^2_b(\R^d)$ be the collection of $C^2$ functions in $\RR^d$
which, along with their partial derivatives of order up to 2, are bounded.
Recall that the infinitesimal generator $\sL$ of the process $X$ is given by
\begin{equation}\label{e:infinitesimal-generator}
    \sL f(x)=\int_{\R^d}\left( f(x+y)-f(x)-y\cdot \nabla f(x) \1_{\{|y|\le\eps\}} \right)\, J(y)dy
\end{equation}
for every $\eps>0$ and $f\in C_b^2(\R^d)$.

\begin{lemma}\label{l:lnew}
There exists $c=c(\phi)>0$ such that for all  $0<r  \le R<\infty$ and
$f\in C^2_b(\R^d)$ with $0\leq f \leq 1$,
$$
\sup_{x \in \R^d}|
\sL f_r(x) |\le   c\, \phi(r^{-2}) \left(
1+\sup_{y}\sum_{j,k} |(\partial^2/\partial
y_j\partial y_k) f(y)| \right) + 2\int_{|z| > R} J(z)dz
$$
 where $f_r(y):=f(y/r)$.
\end{lemma}
\pf
With Lemma \ref{l:integral-estimates-phi} in hand, the proof of this lemma is exactly the same as that of \cite[Lemma 4.2]{KSV7}.
\qed

Similarly, by following the proof of \cite[Lemma 4.10]{KSV7} and using Lemma \ref{l:lnew} instead of \cite[Lemma 4.2]{KSV7}, we obtain the next result.
\begin{lemma}\label{l2.1}
For every $a \in (0, 1)$, there exists $c=c(\phi, a)>0$ such that
for any $r > 0 $ and any open set $D$ with $D\subset B(0, r)$ we have
$$
{\P}_x\left(X_{\tau_D} \in B(0, r)^c\right) \,\le\, c\,
\phi(r^{-2})\ \E_x\tau_D, \qquad x \in D\cap B(0, ar)\, .
$$
\end{lemma}

With the preparation above, we can use  Corollary \ref{c:J-G-a}, Theorem \ref{t:uhi},
Lemma \ref{l:exit-time},  Proposition \ref{p:poisson-kernel-estimate} and  Lemma \ref{l2.1} and
repeat the argument of \cite[Section 5]{KSV7} to get the
following global uniform boundary Harnack principle.
We omit the details here since the proof would be a repetition of the argument in \cite[Section 5]{KSV7}.
Recall that a function $f: \R^d\to [0, \infty)$ is said to be regular harmonic in an open set $U$ with respect to $X$ if
for each $x \in U$,
$f(x)= \E_x\left[f(X(\tau_{U}))\right]$.
\begin{thm}\label{t:ubhp}
Assume {\bf (H1)} and {\bf (H2)}.
There exists $c= c(\phi,d)>0$ such that for every $z_0 \in \R^d$,
every open set $D\subset \R^d$, every $r >0$ and any nonnegative functions $u, v$ in
$\R^d$ which are regular harmonic in $D\cap B(z_0, r)$ with respect to $X$ and
vanish a.e. in $D^c \cap B(z_0, r)$, we have
        $$
        \frac{u(x)}{v(x)}\,\le c\,\frac{u(y)}{v(y)}\, , \qquad \mbox{ for all } x, y\in D\cap B(z_0, r/2).
        $$
\end{thm}

\begin{remark}{\rm
Very recently, the boundary Harnack
principle for (discontinuous) Markov processes (not necessarily L\'evy processes)
on metric measure state spaces is discussed in \cite{BKK}.
In particular in case of a L\'evy processes in $\R^d$,  the boundary Harnack
principle in \cite{BKK} can be stated as follows (see \cite[Theorem 3.5 and Example 5.5]{BKK}):

Let $x_0 \in \R^d$, $0 < r < R $, and let $U \subset B(x_0, R)$ be open. Suppose that $Y$ is
purely discontinuous L\'evy process satisfying \cite[(2.10) and (5.2)]{BKK}.
There exists $c_{(1.1)} = c_{(1.1)}(x_0, r, R)$ such that if $f, g$ are nonnegative functions
on $\R^d$ which are regular harmonic in $U$ with respect to $Y$ and
vanish in $B(x_0, R) \setminus D$,
\begin{equation}\label{e:BKK1}
 f(x)g(y)  \le c_{(1.1)} \, f(y)g(x) \, ,\quad  x, y \in B(x_0, r).
\end{equation}
Condition \cite[(5.2))]{BKK} holds for $X$ by our \eqref{e:J}.
If $d > 2(\delta_2\vee \delta_4)$, \cite[(2.10)]{BKK} holds for $X$ by our \eqref{e:G}.
Comparing with our Theorem \ref{t:ubhp}, the comparison constant
$c_{(1.1)}$ in \eqref{e:BKK1} depends on $x_0$, $r$ and $R$ in general.
It  requires more accurate estimates  to obtain the scale-invariant version of
the boundary Harnack principle, that is,
$c_{(1.1)}$ is independent of $x_0$ and depends on $r$ and  $R$ only through $r/R$.
In fact, in \cite[Example 5.5]{BKK}, it is claimed, without proof, that
one can prove the scale-invariant versions of the boundary Harnack inequalities in \cite{KM2, KSV7} by
checking all dependencies of $c_{(1.1)}$ in \cite[(3.10) and (3.11)]{BKK}.
However, to accomplish this, one needs
to estimate the Green function in order to check Assumption D in \cite{BKK}.
Especially when $X$ is recurrent,  to check Assumption D in \cite{BKK}
one may need upper bounds on the $\alpha$-potential kernel with $\alpha>0$
(see \cite[Proposition 2.3 and the end of the second paragraph of Example 5.5]{BKK}),
which is not discussed in that paper. } \end{remark}


\section{Boundary Harnack principle with explicit decay rate}
Let $D$ be an
open set in $\bR^d$.
For $x\in \bR^d$, let $\delta_{\partial D}(x)$ denote the Euclidean
distance between $x$ and $\partial D$.
Recall that for any
$x\in D$, $\delta_{ D}(x)$ is the Euclidean
distance between $x$ and $D^c$.

In this section we will assume that $D$ satisfies the following types of ball conditions with
radius $R$:
\begin{itemize}
	\item[(i)]
{\it uniform interior ball condition}:
	 for every $x\in D$
with $\delta_{D}(x)< R$ there exists $z_x\in \partial D$ so that
$$
|x-z_x|=\delta_{\partial D
}(x)\ \ \text{ and }\ \
B(x_0, R)\subset D, \quad x_0:=z_x+R\frac{x-z_x}{|x-z_x|};
$$
	\item[(ii)]
{\it uniform exterior ball condition}:
     $D$ is equal to the interior of $\overline{D}$ and for every $y\in \R^d\setminus \overline{D}$
with $\delta_{
\partial D}(y)< R$ there exists $z_y\in \partial D$ so that
$$
|y-z_y|=\delta_{\partial D
}(y)\ \ \text{ and }\ \
B(y_0, R)\subset \R^d\setminus D, \quad y_0:=z_y+R\frac{y-z_y}{|y-z_y|}.
$$
\end{itemize}

The goal of this section is to obtain a global uniform boundary Harnack principle
with explicit decay rate  in open sets in $\bR^d$ satisfying the interior
and exterior ball conditions with radius $R>0$.
This boundary Harnack principle is global in the sense that it holds
for all $R>0$ and the comparison  constant  does not depend on  $R$, and it is uniform
in the sense that it holds for all balls with radii $r \le R$ and
the comparison  constant  depends neither on $D$ nor on $r$.
Throughout the section we assume that {\bf (H1)} and {\bf (H2)} hold.

Let $Z=(Z_t)_{t\ge 0}$ be the one-dimensional subordinate Brownian motion defined
by $Z_t:=W^d(S_t)$. Recall that the potential measure of the ladder height process of
$Z$ is denoted by $V$ and its density by $v$. We also use $V$ to denote the renewal
function of the ladder height process of $Z$. We use the notation
$
\bH:=\{x=(x_1, \dots, x_{d-1}, x_d):=(\tilde{x}, x_d)  \in \bR^d: x_d > 0 \}$ for the half-space.

Define $w(x):=V((x_d)^+)$. Note that $Z_t:=W^d(S_t)$ has a transition density.
Thus, using \cite[Theorem 2]{Sil}, the proof of the next result is the same as
that of \cite[Theorem 4.1]{KSV5}. So we omit the proof.

\begin{thm}\label{t:Sil}
The function $w$ is harmonic in $
\bH$ with respect to $X$ and, for any $r>0$, regular harmonic
in $\R^{d-1}\times (0, r)$ with respect to $X$.
\end{thm}

\begin{prop}\label{c:cforI}
There exists  $c>0$ such that
for all $r>0$ , we have
$$
\sup_{x \in \R^d:\, 0<x_d \le 8 r} \int_{B(x, r)^c \cap
\bH} w(y) j(|x-y|)\, dy \le c\sqrt{\phi(r^{-2})}\, .
$$
\end{prop}
\pf Without loss of generality, we assume $\wt x=0$.
By the substitution $y=x+z$ we see that
$$
\int_{B(x, r)^c \cap \bH} w(y) j(|x-y|)\, dy = \int_{B(0,r)^c\cap \{z_d>-x_d\}} w(z+x)j(z)\, dz=\int_{B(0,r)^c\cap \{z_d>-x_d\}} V(z_d+x_d)j(z)\, dz\, .
$$
The last integral is an increasing function of $x_d$ implying that the supremum is attained for $x_d=8r$.
To conclude, take $x=(\wt{0},8r)$. Then by Theorem 4.1 and \eqref{e:pke-lower},
$$
V(8r)=w(x)=\int_{B(x,r)^c\cap \bH} w(y) K_{B(x,r)}(x,y)\,dy \ge  c_2 \phi(r^{-2})^{-1}\int_{B(x,r)^c\cap \bH} w(y) j(|x-y|)\, dy\, .
$$
Hence,
$$
\int_{B(x,r)^c\cap \bH} w(y) j(|x-y|)\, dy\le  c_3 V(8r)\phi(r^{-2}) \le c_4 \phi(r^{-2})^{1/2}\, .
$$
\qed

For a function $f:\R^d\to  \R$ and $x\in \R^d$ we define
$$
     \sA  f(x)
     :=\lim_{\eps \downarrow 0} \int_{\{y\in \bR^d: |x-y| > \eps\}} \left(f(y)-f(x)\right)j(|y-x|)\, dy\, ,
$$
and use  $\mathfrak{D}_x(\sA)$ to denote the family of all functions $f$ such that
$\sA  f(x)$
exists and is finite.
It is well known that $C^2_c(\RR^d)\subset \mathfrak{D}_x(\sA)$ for every $x \in \RR^d$
and that, by the rotational symmetry of $X$, $\sA$ restricted to $C^2_c(\RR^d)$
coincides with the infinitesimal generator $\sL$ of $X$ which is given in \eqref{e:infinitesimal-generator}.

Using \cite[Corollary 13.3.8]{KSV3}, Theorem \ref{t:Sil},  \eqref{e:j-decay} and
\eqref{e:doubling-condition}, the proof of the next result is the same as these of
\cite[Proposition 4.3]{KSV5} and \cite[Theorem 3.4]{KSV4}, so we omit the proof.

\begin{thm}\label{c:Aw=0}
For any $x\in
\bH$, $w\in \mathfrak{D}_x(\sA)$ and $\sA w(x)=0$.
\end{thm}

Before we prove our main technical lemma, we first do some preparations.

\begin{lemma}\label{4.3.5}
If $f, g: (0, \infty)\to (0,\infty)$ are non-increasing, then for any $M>0$ and any $x:[0,M]\to \R$ we have
$$
    \int_0^M \int_0^M f(s)g(r+|s-x(r)|)\, dr\, ds \le 2\int_0^{3M/2} F(u) g(u)\, du\, ,
    $$
    where $F(u)=\int_0^u f(s)\, ds$.
\end{lemma}

\pf
Without loss of generality we may assume that $g$ is right continuous. Then the inverse $g^{-1}(\lambda):=\sup\{x:\, g(x)\ge \lambda\}$ has the property that $g(x)\ge \lambda$ if and only if $x\le g^{-1}(\lambda)$.
    Let $h(s):=g(r+|s-x(r)|)$, $s\in [0,M]$. Then
    \begin{eqnarray*}\lefteqn{\big|\{s\in [0,M]:\, h(s)>\lambda\}\big| = \big|\{s\in [0,M]:\, g(r+|s-x(r)|)>\lambda\}\big|}\\
    &=&\big|\{s\in [0,M]:\, r+|s-x(r)|\le g^{-1}(\lambda)\}\big| \\
    &=&\big|\{s\in [0,M]:\, |s-x(r)|\le (g^{-1}(\lambda)-r)^+\}\big| \\
    &\le & 2(g^{-1}(\lambda)-r)^+\ .
    \end{eqnarray*}
Hence, the rearrangement $\{s\in [0,M]:\, h(s)>\lambda\}^*$ is contained in $[0,2(g^{-1}(\lambda)-r)^+]$.
    Further note that $s\le 2(g^{-1}(\lambda)-r)^+$ is equivalent to $r+\frac{s}{2}\le g^{-1}(\lambda) $,
     which in turn is equivalent to $g(r+s/2)\ge \lambda$. Therefore the non-increasing rearrangement of $h$ satisfies
    \begin{eqnarray*}
    h^*(s)&=&\int_0^{\infty} \1_{\{h>\lambda\}^*}(s)\, d\lambda \le \int_0^{\infty} \1_{[0,2(g^{-1}(\lambda)-r)^+]}(s)\, d\lambda 
    =\int_0^{\infty}
    \1_{[0,g(r+s/2)]}(\lambda)\, d\lambda\\
    & =&\int_0^{g(r+s/2)}\, d\lambda \le g(r+s/2)\, .
    \end{eqnarray*}
Therefore, by the rearrangement inequality (see \cite[Chapter 3]{LL}),

   \begin{eqnarray*}
   \int_0^M f(s)g(r+|s-x(r)|)\, ds =\int_0^M f(s)h(s)\, ds \le \int_0^M f(s)h^*(s)\, ds\le \int_0^M f(s)g(r+s/2)\, ds\, .
   \end{eqnarray*}
   Finally,
   \begin{eqnarray*}
   \lefteqn{\int_0^M \int_0^M f(s)g(r+|s-x(r)|)\, dr\, ds \le \int_0^M \int_0^M f(s)g(r+s/2)\, dr\, ds}\\
   &\le& \int_0^{3M}\int_{s/2}^{3M/2} f(s)g(u)\, du\, ds =\int_0^{3M/2}\left(\int_0^{2u}f(s) \,ds \right)g(u)\, du \\
   &=&\int_0^{3M/2}F(2u)g(u)\, du\le 2\int_0^{3M/2} F(u)g(u)\, du\, .
   \end{eqnarray*}

\qed

\begin{lemma}\label{l:main0}
Let $D$ be an open set in $\bR^d$ satisfying
the interior and exterior ball conditions with
radius $1$. Fix $x \in D$ with $\delta_D(x) <1/8$ and
let $x_0 \in\partial {D}$ be such that $\delta_{D}(x)=|x-x_0|$ and  $CS_{x_0}$ be a coordinate system such that $x=( \wt 0, x_d)$ and $x_d >0$.
There exists $c>0$  independent of $D$ and $x$ such that
for  every positive non-increasing functions $\nu$ and $\vartheta$ on $(0, \infty)$ and $\varTheta (r)=\int_0^r \vartheta (s)ds$
\begin{align}
\label{e:ndfe1}
\int_{B(x, 1/8)}|
\varTheta (\delta_D(z))-\varTheta (\delta_{H^+} (z))|\, \frac{\nu(|z-x|)}{|z-x|^d}\, dz \le c \int_{0}^1  \varTheta (2r) \nu(r)\, dr\ ,
\end{align}
where $H^+:=\left\{z=(\wt z, \, z_d) \mbox{ in } CS_{x_0} :z_d>0 \right\}$.
\end{lemma}
\pf
In this proof we assume that $d \ge 2$, the case $d=1$ being simpler.
By the interior and exterior ball conditions with radius $1$,
$$
\{ z=(\wt z, \, z_d) \in B(0, 1/2):  z_d  \ge \psi (|\wt z|) \}\subset B(0, 1/2)\cap {D}\subset
\{ z=(\wt z, \, z_d) \in B(0, 1/2):  z_d  \ge -\psi (|\wt z|)  \},
$$
where
$ \psi(r):=1-\sqrt{ 1-r^2}$.

Define
\begin{eqnarray*}
     A&:=&\{z=(\widetilde{z},z_d) \in ({D}\cup H^+) \cap B(x, {1}/{8}):
    -\psi (|\wt z|)   \le z_d  <\psi (|\wt z|) \}\, ,\\
    F&:=&\{z\in B(x,\tfrac{1}{8})\colon z_d>\psi (|\wt z|) \} \,.
\end{eqnarray*}
Then
\begin{align*}
&\int_{B(x, 1/8)}|
\varTheta (\delta_D(z))-\varTheta (\delta_{H^+} (z))| \frac{\nu(|z-x|)}{|z-x|^d} dz
\\
 \le &\int_{A}
\varTheta (\delta_D(z))+\varTheta (\delta_{H^+} (z)) \frac{\nu(|z-x|)}{|z-x|^d} dz
+\int_{F} |
\varTheta (\delta_D(z))-\varTheta (\delta_{H^+} (z))| \frac{\nu(|z-x|)}{|z-x|^d} dz =:I+II.
\end{align*}
Let $E=B((\wt 0, -1), 1)^c$. Then
\begin{align*}
    I
    \leq 2\int_0^\frac{1}{8}\int_{|\wt z|=r}\mathbf{1}_{\{z=(\widetilde{z},z_d) : |\wt z|=r, -\psi(r)\le z_d <\psi(r)  \}}(z)\varTheta\left(\delta_{E} (z)\right)
    \frac{\nu(\sqrt{r^2+|z_d-x_d|^2})}{(r^2+|z_d-x_d|^2)^{d/2}}\,m_{d-1} (dz)\, dr\,,
\end{align*}
where $m_{d-1}$ is the surface measure, that is, the $(d-1)$-dimensional Lebesgue measure.
Noting that $1-\sqrt{1-|\wt z|^2}\leq {|\wt z|^2}{}=r^2$ for $|\wt z|=r$, we obtain
$$
    m_{d-1}\big(\big\{z=(\widetilde{z},z_d) : |\wt z|=r, -1 +\sqrt{ 1-r^2}\le z_d <1 -\sqrt{ 1-r^2} \big\}\big)\leq c_1 {r^d}\ \ \text{ for }\ \ r\leq \tfrac{1}{8}\, .
$$
Since $\varTheta$ is increasing and $1-\sqrt{1-|\wt z|^2}\leq {|\wt z|^2}{}\leq {|\wt z|}{}$ we deduce
$\varTheta\left(\delta_{E} (z)\right)\le \varTheta\left(2 \psi (|\wt z|) \right) \leq \varTheta(2|\wt z|)$.
By using that $\nu$ is decreasing we get
\begin{align*}
  I  & \leq \  \int_0^\frac{1}{4}\int_{|\wt z|=r}\mathbf{1}_{\{z\colon |\wt z|=r,\ -1 +\sqrt{   1-r^2}\le z_d < 1 -\sqrt{ 1-r^2}\}}(z)
  \varTheta(2r)\nu(r)r^{-d}\,m_{d-1}(dz)\,dr \\
   &\leq c_2 \int_0^\frac{1}{4}  \varTheta(2r)\nu(r) dr.
\end{align*}

In order to estimate $II$, we consider two cases. First,
if $0 <z_d=\delta_{_{H^+}}({z}) \le  \delta_{ D}({z})$, then using the exterior ball condition
and the fact that $\delta_{ D}({z})$ is smaller than the vertical distance from $z$ to the exterior
of the ball which is equal to $z_d+1-\sqrt{1-|\wt{z}|^2}\le z_d +|\wt{z}|^2$, we get that
\begin{align}
\varTheta(\delta_{ D} (z))-\varTheta(\delta_{H^+}(z)) \le  \varTheta(z_d+|\wt z|^2) -\varTheta(z_d) =
\int_{z_d}^{z_d+|\wt z|^2} \vartheta(t)dt \le |\wt z|^2
\vartheta(z_d),  \label{e:KG}
\end{align}
since $\vartheta$ is decreasing.

If $z_d=\delta_{_{H^+}}({z}) >  \delta_{ D}({z})$ and $z \in F$,
using the fact that
$\delta_{ D}({z})$ is greater than or equal to the distance between $z$
and the graph of $ \psi$ (by the interior ball condition) and
 \begin{align*}
z_d-1+\sqrt{ |\wt z|^2+(1-z_d)^2} &= \tfrac{|\wt z|^2}
{\sqrt{ |\wt z|^2+(1-z_d)^2} + (1-z_d)}\,\le\, \tfrac{ |\wt z|^2} {2 (1-z_d)}
\le
 |\wt z|^2,\quad \forall z \in F,
\end{align*}
we obtain
 \begin{align}
\varTheta(\delta_{H^+}(z))-\varTheta\left(\delta_{ D} (z)\right)
\le\int^{z_d}_{1-\sqrt{ |\wt z|^2+(1-z_d)^2}} \vartheta(t)dt\le    |\wt
z|^2
\,\vartheta\left(1-\sqrt{ |\wt z|^2+(1-z_d)^2}\right). \label{e:KG2}
\end{align}
By \eqref{e:KG} and \eqref{e:KG2},
\begin{align*}
II
\,\le\, \int_{F}
|\wt z|^2 \left(\vartheta(z_d) \vee \vartheta(1-\sqrt{ |\wt z|^2+(1-z_d)^2})\right) \frac{\nu(|z-x|)}{|z-x|^d} dz\,
.
\end{align*}
Since
\[
F\subset  \{z=(\wt z, z_d)\in \R^d: \ |\widetilde{z}|< {1}/{8}
\hbox{ and }  \psi (r) < z_d \leq {1}/{4}\},
\]
switching to polar coordinates for $\wt z$ and reversing the order of integration   we get
\begin{align*}
II\leq c_{3} \int_{0}^{\frac{1}{8}} \int_{\psi (r)}^{ \psi (r)+1/4}
( \vartheta(z_d) \vee \vartheta(1-\sqrt{ r^2+(1-z_d)^2}))
\frac{r^{d}\nu(\sqrt{r^2+|z_d-x_d|^2})}{(r^2+|z_d-x_d|^2)^{d/2}} dz_d dr\ .
\end{align*}
Writing $s=z_d-\psi (r)$ gives
\begin{align*}
II\leq c_{3}
\int_{0}^{\frac{1}{8}}
\int_{0}^{ 1/4}
( \vartheta(\psi (r) +s) \vee  \vartheta(1-\sqrt{
 1-2s \sqrt{1-r^2} +s^2
 })) r^{d} \frac{\nu(\sqrt{r^2+|s+\psi (r) -x_d|^2})}{(r^2+|s+\psi (r)-x_d|^2)^{d/2}} ds dr.
\end{align*}
Since
$$
1-\sqrt{
 1-2s \sqrt{1-r^2} +s^2
 } \ge (2s \sqrt{1-r^2} -s^2)/2=(2 \sqrt{1-r^2} -s)s/2 \ge s/2
 $$
 and
 $\sqrt{r^2+|z_d-x_d|^2} \ge (r+|z_d-x_d|)/2$, we have
 \begin{align*}
II\leq c_{3}
\int_{0}^{\frac{1}{8}}
\int_{0}^{ 1/4}
 \vartheta(s/2)
  \nu((r+|s+\psi (r) -x_d|)/2) ds dr.
\end{align*}
Thus by Lemma \ref{4.3.5},
$$ II\le 2c_3
\int_{0}^{3/8}  \varTheta (r/2) \nu(r/2)dr \le c_4
\int_{0}^{1}  \varTheta (r) \nu(r)dr .$$
\qed

\begin{prop}\label{l:main}
Let $D$ be an open set in $\bR^d$ satisfying
the interior and exterior ball conditions with
radius $R$. Fix $Q \in \partial D$  and define
$$
h(y):=V\left(\delta_{D} (y)\right){\bf 1}_{D\cap B(Q, R/2)}(y).
$$
There exists $C_1= C_1(\phi)>0$  independent of $Q$  and $R$
(and $D$) such that
$h\in \mathfrak{D}_x(\sA)$ for every $x\in {D}\cap B(Q, R/8)$
and
\begin{equation}\label{e:main-lemma}
| \sA h(x)| \le C_1  \sqrt{\phi(R^{-2})} \quad \text{ for all }
x \in {D}\cap B(Q, R/8)\, .
\end{equation}
\end{prop}
\pf
In this proof we assume that $d \ge 2$, the case $d=1$ being simpler.
We fix $x \in {D}\cap B(Q, {R}/{8})$ and
let $x_0 \in\partial {D}$ be such that $\delta_{D}(x)=|x-x_0|$.
We may assume, without loss of generality, that $x_0=0$, $x=( \wt 0, x_d)$ and $x_d>0$.
Note that, since $|y-Q| \le |y-x| +|x-Q| \le R/4$ for $y \in B(x, {R}/{8})$, we have
    \begin{equation}\label{e:BE}
B(x, {R}/{8}) \cap {D} \subset B(Q,   {R}/{4})\cap {D} \, .
\end{equation}

Let $h_{x}(y):=V(\delta_{_{\bH}}(y)).$
Note that $h_{x}(x)=h(x)$.
Since $\delta_{_{\bH}}(y)=(y_d)^+$,
it follows from Theorem \ref{c:Aw=0} that $\sA h_x$ is well defined in $\bH$ and
\begin{equation}\label{e:hz}
    \sA  h_{x}(y)=0, \quad  \forall y\in \bH.
\end{equation}
We show now that $\sA (h-h_x)(x)$ is well defined. For each small $\varepsilon >0$ we have that
\begin{align*}
    &\int_{\{y \in {D} : \, |y-x|>\varepsilon\}}
    |h(y)-h_{x}(y)|j(|y-x|)\ dy \nonumber\\
    \leq & \int_{B(x,\frac{R}{8})^c} (h(y)+h_{x}(y)) j(|y-x|)dy
   +\int_{B(x,\frac{R}{8})} {|h(y)-h_{x}(y)|}j(|y-x|) dy=:I_1+I_2.
\end{align*}

We claim that
\begin{equation}\label{e:I1234}
I_1+I_2\le C_1\sqrt{\phi(R^{-2})}
\end{equation}
for some positive constant $C_1>0$.
Since \eqref{e:I1234} implies that
\begin{align*}
    &{\bf 1}_{\{y \in {D}\cup \bH: \, |y-x|>\varepsilon\}}{|h(y)-h_{x}(y)|}j(|y-x|) \\
     &\le {\bf 1}_{B(x,\frac{R}{8})^c} (h(y)+h_{x}(y)) j(|y-x|)+{\bf 1}_{B(x,\frac{R}{8})}{|h(y)-h_{x}(y)|}j(|y-x|)  \in L^1(\R^d)\, ,
\end{align*}
by the dominated convergence theorem the limit
$$
    \lim_{\varepsilon\downarrow 0}\int_{\{y \in {D} \cup \bH: |y-x|>\varepsilon\}}{(h(y)-h_{x}(y))}j(|y-x|)\, dy
$$
exists, and hence $\sA(h-h_x)(x)$ is well defined and $|\sA(h-h_x)(x)|\le C_1\sqrt{\phi(R^{-2})}$. By linearity and \eqref{e:hz}, we get that $\sA h(x)$ is well defined and $|\sA h(x)|\le C_1\sqrt{\phi(R^{-2})}$. Therefore, it remains to prove \eqref{e:I1234}.

Since $h(y)=0$ for $y \in B(Q, R)^c$, it follows that
\begin{align}
    I_1 \le&
    \int_{B(x,\frac{R}{8})^c} V(y_d) j(|y-x|)dy+  V(R)\int_{B(x,\frac{R}{8})^c}j(|y-x|)dy \nn\\
   \le &\sup_{z \in \R^d:\ 0<z_d <R} \int_{B(z,  \frac{R}{8})^c \cap \bH}    V(y_d) j(|z-y|) dy
   + V(R)\int_{B(0,  \frac{R}{8})^c}j (|y|)dy\nn\\
    =:&K_1+V(R) K_2. \nonumber
\end{align}
By Proposition \ref{c:cforI} we have that $K_1 \le c \sqrt{\phi(R^{-2})}$.  Moreover,
by Lemma \ref{l:j-g-upper},
  \eqref{e:V-phi-a} and \eqref{e:ie-2},
  $$
  V(R)\int_{B(0,  \frac{R}{8})^c}j (|y|)dy \le  c_2 V(R)\int_{R/8}^\infty  r^{-1} \phi(r^{-2})dr \le c_3 V(R)\phi(R^{-2}/64)  \le c_4 \sqrt{\phi(R^{-2})}.
  $$

For $I_2$, we use scaling. Let $x^R=R^{-1} x$ and $ \wh D:=\{z: Rz \in D \}$.
Then  by \eqref{e:jumping-function-a} and \eqref{e:Va} we have
\begin{align*}
  I_2 &=\int_{\{y\in B(x,\tfrac{R}{8})\colon y_d>R -\sqrt{ R^2-|\wt y|^2}\}}
   {|V(\delta_{D} (y))-V(\delta_{\bH}(y))|}j(|y-x|) dy\\
          &=\sqrt{\phi(R^{-2})} \int_{B(x^R,1/8)} {|V^R(\delta_{ \wh D} (z))-
         V^R(\delta_{\bH}(z))|}j^R(|z-x^R|) dz=: \sqrt{\phi(R^{-2})}\,
         \wh I_2\ .
\end{align*}

Using \eqref{e:J-a} and \eqref{e:V-phi-a},
$$
\wh I_2  \le c_5 \int_{B(x^R,1/8)} |V^R(\delta_{ \wh D} (z))-
         V^R(\delta_{\bH}(z))|
          \frac{       \phi^R(|z-x^R|^{-2})}{   |z-x^R|^d}             dz\ .
$$
       Finally by   Lemma \ref{l:phi-property}, Lemma \ref{l:main0}, \eqref{e:V-phi-a} and \eqref{e:sc1},
       $$
    \wh I_2  \le c_6    \int_0^1  V^R(2r)\phi^R(r^{-2})
    \le c_7 \int_0^1 \sqrt{\frac{\phi(R^{-2}(2r)^{-2})}{\phi(R^{-2})}}dr \le c_8\int^1_0r^{-(\delta_2\vee\delta_4)}dr < \infty.
    $$
\qed

\begin{thm}\label{L:2}
\noindent
(a) There exist
$a=a(\phi)\in (0,1)$
 and $c_1=c_1(\phi)>0$ such that for every open set $D$
satisfying the interior and exterior ball conditions with
radius $R>0$, any $r \le aR$ and $Q \in \partial D$,
\begin{equation}\label{e:L:3}
\E_x\left[ \tau_{D \cap B(Q, r)}\right] \le c_1 V(r)V(\delta_{D} (x)),\qquad
\hbox{for every } x \in  D \cap B(Q, r)   \, .
\end{equation}

\noindent
(b)
There exists $c_2=c_2(\phi)>0$  such that for every open set $D$
satisfying the interior and exterior ball conditions with
radius $R>0$, $r \in (0, R]$,
$Q\in \partial D$ and any nonnegative function $u$ in $\R^d$
which is harmonic in $D \cap B(Q, r)$ with respect to $X$ and
vanishes continuously on $ D^c \cap B(Q, r)$, we have
\begin{equation}\label{e:bhp_m}
\frac{u(x)}{u(y)}\,\le c_2\,
\sqrt{\frac{\phi(\delta_{D}(y)
^{-2})}{\phi(\delta_{D}(x)
^{-2})}}\ ,
\qquad \hbox{for every } x, y\in  D \cap B(Q, \tfrac{r}{2}).
\end{equation}
\end{thm}
\pf
Without loss of generality, we assume $Q=0$.
Define
$$
h(y) := V(\delta_{D}(y)) {\bf 1}_{B(0, R/2) \cap D}(y)\, .
$$
Let $f$ be a non-negative smooth radial function such that $f(y)=0$ for $|y|>1$ and
$\int_{\R^d} f (y) dy=1$. For $k\geq 1$, define $f_k(y)=2^{kd} f (2^k y)$ and
$$
h^{(k)}(z):= ( f_k*h)(z) :=\int_{\R^d}f_k (y) h(z-y)dy\, ,
$$
and for
$\lambda \ge 8$ let $B^\lambda_k:=\left\{y \in D \cap B(0, \lambda^{-1}R): \delta_{D \cap B(0, \lambda^{-1}R)    }(y) \ge 2^{-k}\right\}$.
Since $h^{(k)}$ is a $C^{\infty}$ function, $\sA h^{(k)}$ is well defined everywhere.
Then by the same argument as that in \cite[Lemma 4.5]{KSV5}, we have for large $k$
\begin{equation} \label{e:*333}
-C_1 \sqrt{\phi(R^{-2})}\le \sA h^{(k)} \leq C_1 \sqrt{\phi(R^{-2})}\quad \text{ on } B^\lambda_k\, ,
\end{equation}
where $C_1$ is the constant from Proposition \ref{l:main}.

Since $h^{(k)}$ is in $C^\infty_c(\R^d)$ and that $\sA$ restricted to $C^\infty_c$ coincides with the infinitesimal generator $\sL$ of
the process $X$, by Dynkin's formula, with $\sigma(\lambda,k):= \tau_{B^\lambda_k}$
\begin{equation} \label{e:*334}
\E_x\int_0^{\sigma(\lambda,k)}  \sA h^{(k)}(X_t) dt =\E_x[h^{(k)}(X_{\sigma(\lambda,k)})]
-h^{(k)}(x)\, .
\end{equation}
Using \eqref{e:*333}--\eqref{e:*334} and then letting $k\to \infty$ we obtain
that for all $\lambda \ge 8$ and $x \in D \cap B(0, \lambda^{-1}R)$,
\begin{align}
    &V(\delta_{D}(x)) \,=\, h(x) \ge \E_{x}\left[h\left(X_{ \tau_{D \cap B(0, \lambda^{-1}R)}} \right)  \right] - C_1 \sqrt{\phi(R^{-2})}\E_x\left[\tau_{D \cap B(0, \lambda^{-1}R)}\right] \label{e:ss}
\end{align}
and
\begin{align}
    &V(\delta_{D}(x)) - C_1 \sqrt{\phi(R^{-2})}\, \E_x\left[\tau_{D \cap B(0, \lambda^{-1}R)}\right] \le \E_{x}\left[h\left(X_{ \tau_{D \cap B(0, \lambda^{-1}R)}} \right)  \right]  \label{e:ss1}.
\end{align}
Since
$$
j(|y-z| ) \ge j(|y|+|z|) \ge   j(2|y|) \ge  c_1 j(|y|),  \quad \forall (z,y)  \in (D \cap B(0, \lambda^{-1}R)) \times B(0,  \lambda^{-1} R)^c,
$$
we get
\begin{align}
    &\int_{ (B(0,R)\setminus B(0, \lambda^{-1}R) )\cap D} \int_{D \cap B(0, \lambda^{-1}R)} G_{D \cap B(0, \lambda^{-1}R)}(x,z) j(|z-y|)dz    V(\delta_{D}(y)) dy \nonumber \\
    &\ge  c_1  \E_x\left[\tau_{D \cap B(0, \lambda^{-1}R)}\right] \int_{ (B(0,R)\setminus B(0, \lambda^{-1}R) )\cap D}    j(|y|)V(\delta_{D}(y)) dy \label{e:new23}.
\end{align}
The remainder of the proof is written for $d\ge 2$. The interpretation in the case $d=1$ is obvious.

By the interior ball condition with radius $R$, we may assume, without loss of generality, that
$$
B(0, R)\cap {D}\supset
\{ y=(\wt y, \, y_d) \in B(0, R)
: R -\sqrt{ R^2-|\wt y|^2}< y_d \}.
$$
For $ y \in B(0, R)$ with $2|\wt y| <  y_d$ we have
\begin{align}
\delta_{D}(y) \ge
R-\sqrt{ |\wt y|^2+(R-y_d)^2}
\ge R-\sqrt{R^2-2Ry_d+(5/4)y^2_d}\nn\\
\ge (2Ry_d-(5/4)y^2_d)/(2R)=y_d(1-(5/8)y_d/R)
\ge\tfrac{3y_d}{8 }\ge \tfrac{3|y|}{4\sqrt 5 }  \, .
\label{e:new14}
\end{align}
Thus, by changing into polar coordinates and using \eqref{e:new14}, we have
\begin{align*}
  \int_{ (B(0,R)\setminus B(0, \lambda^{-1}R) )\cap D}    j(|y|)V(\delta_{D}(y)) dy
    & \ge \int_{ \{(\wt y, y_d): 2 |\wt y|< y_d,   \lambda^{-1} R<|y| <R\}} j(|y|)V(\tfrac{3|y|}{4 \sqrt 5 }  ) dy\\
    &\ge c_2 \int_{  \lambda^{-1} R}^R    j(r)V(\tfrac{3r}{4 \sqrt 5 }  )  r^{d-1}\, dr\, .
\end{align*}
By Lemma \ref{l:phi-property}, \eqref{e:Berall1}, \eqref{e:J} and \eqref{e:V-phi-a}, we have that
for $r <R$,
\begin{align*}
j(r)V(\tfrac{3r}{4 \sqrt 5 })  r^{d-1} \ge c_3 r^{-3} \frac{\phi'(r^{-2})} { \phi(r^{-2})^{1/2}}
=c_4 \frac{d}{dr}(-\phi(r^{-2})^{1/2})\, .
\end{align*}
Thus
\begin{align}
   \int_{ (B(0,R)\setminus B(0, \lambda^{-1}R) )\cap D}    j(|y|)V(\delta_{D}(y)) dy
    \ge c_4 (\phi( R^{-2} \lambda^{2})^{1/2} - \phi( R^{-2} )^{1/2})\, . \label{e:new25}
\end{align}
Now combining \eqref{e:new23} and \eqref{e:new25}, we get
\begin{align}
    \E_{x}\left[h\left(X_{ \tau_{D \cap B(0, \lambda^{-1}R)}} \right) \right] \ge
    c_5 \left(\phi( R^{-2} \lambda^{2})^{1/2} - \phi( R^{-2} )^{1/2}\right)\E_x\left[\tau_{D \cap B(0, \lambda^{-1}R)}\right]\, .\label{e:sdf311}
\end{align}
Thus by \eqref{e:ss} for every $x \in D \cap B(0, \lambda^{-1}R)$,
\begin{align}
V(\delta_{D}(x))\ge \left(c_5\phi( R^{-2} \lambda^{2})^{1/2} -
(c_5+C_1) \phi( R^{-2} )^{1/2}\right) \E_x[\tau_{D \cap B(0, \lambda^{-1}R)}]\, . \label{e:ss4}
\end{align}

Without loss of generality we assume $a_5 <1$ (the constant in \eqref{e:sc1}). Let
$
\lambda_0 :=(2 a_5^{-1/2}(1+C_1/c_5) )^{1/(\delta_1 \wedge  \delta_3)} \vee 8.
$
Recall from \eqref{e:sc1} that
\begin{align}
\phi( t) \ge a_5 s^{\delta_1 \wedge  \delta_3} \phi(s^{-1} t) \quad \text{ for every } s \ge 1 \quad \text { and } t>0. \label{e:ss711}
\end{align}
 Applying this with $t= \lambda_0^2 R^{-2}$ and $s= \lambda_0^2 \ge1$,  we get that for $\lambda\ge \lambda_0$,
$$
    \phi( R^{-2} \lambda^{2}) \ge \phi( R^{-2} \lambda_0^{2}) \ge a_5(\lambda_0^2)^{\delta_1 \wedge  \delta_3}
    \phi( R^{-2} ) \ge   4 (1+C_1/c_5)^2 \phi( R^{-2} )\, .
$$
Hence for every $\lambda \ge \lambda_0$
\begin{align}
    & c_5\phi( R^{-2} \lambda^{2})^{1/2} -(c_5+C_1) \phi( R^{-2} )^{1/2} \ge
    \frac{c_5}{2} \phi( R^{-2} \lambda^{2})^{1/2} \label{e:ss6}.
\end{align}
Combining \eqref{e:ss4} and \eqref{e:ss6}, we have proved part (a) of the theorem with $a=\lambda_0^{-1}$.

To prove  (b), we first consider estimates on $h$ first.
Combining \eqref{e:ss1} and \eqref{e:ss4}, we get
\begin{align}
& \E_{x}\left[h\left(X_{ \tau_{D \cap B(0, \lambda^{-1}R)}} \right)  \right] \ge    V(\delta_{D}(x)) - C_1 \sqrt{\phi(R^{-2})}\E_x\left[\tau_{D \cap B(0, \lambda^{-1}R)}\right]\nonumber\\
&\ge    V(\delta_{D}(x)) \left(1- C_1\tfrac{ \sqrt{\phi(R^{-2})}}{c_5\phi( R^{-2} \lambda^{2})^{1/2} -
(c_5+C_1) \phi( R^{-2} )^{1/2} } \right).\label{e:ss41}
\end{align}

Let
$$
\lambda_1 :=\big((3C_1 +c_5)/(
c_5a_5^{1/2})\big)^{{1}/{(\delta_1 \wedge  \delta_3)}} \vee \lambda_0.
$$
 Applying \eqref{e:ss711} with $t= \lambda_1^2 R^{-2}$ and $s= \lambda_1^2 \ge1$,  we get that for $\lambda\ge \lambda_1$,
$$
    \phi( R^{-2} \lambda^{2}) \ge \phi( R^{-2} \lambda_1^{2}) \ge a_5(\lambda_1^2)^{\delta_1 \wedge  \delta_3} \phi( R^{-2} ) =  ( 3C_1 +c_5)^2c_5^{-2}\phi( R^{-2} )\, .
$$
Hence for every $\lambda \ge \lambda_1$,
\begin{align}
    &  2C_1\sqrt{\phi(R^{-2})} \le c_5\phi( R^{-2} \lambda^{2})^{1/2} -
    (c_5+C_1) \phi( R^{-2} )^{1/2} \label{e:ss61}.
\end{align}
Combining \eqref{e:ss41}--\eqref{e:ss61}, we have for every $x \in D \cap B(0, \lambda^{-1}R)$,
\begin{align}
\E_{x}\left[h\left(X_{ \tau_{D \cap B(0, \lambda^{-1}R)}} \right)  \right] \ge  \frac12 V(\delta_{D}(x)). \label{e:ss71}
\end{align}
Moreover, by \eqref{e:ss}, \eqref{e:ss4} and \eqref{e:ss6}, for every $\lambda \ge \lambda_0$ and $x \in D \cap B(0, \lambda^{-1}R)$,
\begin{align}
\E_{x}\left[h\left(X_{ \tau_{D \cap B(0, \lambda^{-1}R)}} \right)  \right] &\le V(\delta_{D}(x))+  C_1 \sqrt{\phi(R^{-2})}\E_x\left[\tau_{D \cap B(0, \lambda^{-1}R)}\right]\nonumber\\
&\le V(\delta_{D}(x))\left(1+  \tfrac{2C_1}{c_5} \sqrt{\tfrac{\phi(R^{-2})}
   {\phi( R^{-2} \lambda^{2})}  }\right) \le  (1+  {2C_1}/{c_5})  V(\delta_{D}(x))
. \label{e:ss81}
\end{align}

Now we assume
$r \in (0, R]$.  Let $u$ be a nonnegative function in $\R^d$
which is harmonic in $D \cap B(0, r)$ with respect to $X$ and
vanishes continuously on $ D^c \cap B(0, r)$.
Note that
$0<r /(2\lambda_1) <r/\lambda_1 =R (R\lambda_1/r)^{-1}$.
Thus by applying Theorem \ref{t:ubhp}  to $u$ and $v(x):=\E_{x}[h(X_{ \tau_{D \cap B(0, \lambda_0^{-1}r)}} )  ]$ first  and then by applying \eqref{e:ss71}--\eqref{e:ss81} (with $\lambda=R\lambda_1/r$),  we obtain that
for every $x, y\in  D \cap B(0, {r}/(2\lambda_1))$,
$$
\frac{u(x)}{u(y)}\,\le c_6\, \frac{v(x)}{v(y)}
\,\le c_7\, \frac{V(\delta_{D}(x))}{
V(\delta_{D}(y))}
 \,\le\, c_8\,
\sqrt{\frac{\phi(\delta_{D}(y)
^{-2})}{\phi(\delta_{D}(x)
^{-2})}}.
$$
When $x$ or $y$ in   $D \cap (B(0, {r}/2)\setminus B(0, {r}/(2\lambda_1)))$,
we first use the standard chain argument and then apply the above result.
\qed


\section{Heat kernel estimates in the half-space}

Recall that $p(t, x, y)$ is the transition density of $X$ and $\Phi$ stands for the function
$\Phi(r)=1/\phi(r^{-2})$, $r>0$.
We use $\Phi^{-1}(r)$ to denote the inverse function of $\Phi$.
Since $X$ satisfies \cite[(1.4), (1.13) and (1.14)]{CK},
by \cite[Theorem 1.2]{CK} the following estimates for $p(t, x, y)$ are valid:
there exists $c_1>0$ such that for all $(t, x, y)\in (0, \infty)\times \R^d\times \R^d$,
\begin{align}\label{stssbound}
c^{-1}_1\left(\frac1{(\Phi^{-1}(t))^d}\wedge tJ(x, y)\right)
\le p(t, x, y) \le c_1 \left(\frac1{(\Phi^{-1}(t))^d}\wedge tJ(x, y)\right).
\end{align}

It is known (see \cite{CK}) that the killed process $X^D$ has a
transition density $p_D(t, x, y)$ with respect to the Lebesgue
measure that is jointly H\"older continuous.
In a recent preprint \cite{CKS6}, sharp two-sided estimates on $p_D(t, x, y)$ for
bounded open sets have been established for subordinate Brownian motions
under weaker conditions.

The goal of this section is to get sharp two-sided estimates for
$p_{\bH}(t, x, y)$, and, as a consequence,
sharp two-sided estimates of the Green function $G_{\bH}(x,y)$.

\begin{lemma}\label{l:u_near}
There exists $c=c( \phi )>1$ such that for every  $(t,x,y) \in (0,, \infty) \times \bH
\times \bH$,
$$
 p_{\bH}(t,x,y) \,\le\, c\, (\Phi^{-1}(t))^{-d}
\left(\sqrt{\frac{\Phi(\delta_{\bH}(x))}{t}} \wedge 1\right)\left(\sqrt{\frac{\Phi(\delta_{\bH}(y))}{t}} \wedge 1\right)  .
$$
\end{lemma}

\pf
Let
$c(t):=\sup_{z,w \in  \R^d} p(t/3,z,w)$.
By the semigroup property and symmetry,
$$
 p_{\bH}(t,x,y) =\int_{\bH} \int_{\bH} p_{\bH}(t/3,x,z)
p_{\bH}(t/3,z,w) p_{\bH}(t/3,w,y) dzdw
 \le c (t) \, \P_{x}(\tau_{\bH}>t/3)
\P_{y}(\tau_{\bH}>t/3).
$$
Now the lemma follows from
Lemma \ref{l:phi-property}, \eqref{stssbound} and \cite[Theorem 4.6]{KMR}.
\qed

The next lemma and its proof are given in \cite{CKS5} (also see
\cite[Lemma 2]{BGR} and \cite[Lemma 2.2]{CKS}).

\begin{lemma}\label{l:gen1}
Suppose that $U_1,U_3, E$ are open subsets of $\bR^d$ with $U_1,
U_3\subset E$ and ${\rm dist}(U_1,U_3)>0$. Let $U_2 :=E\setminus
(U_1\cup U_3)$.  If $x\in U_1$ and $y \in U_3$, then for all $t >0$,
\begin{equation}\label{eq:ub}
p_{E}(t, x, y) \le \P_x\left(X_{\tau_{U_1}}\in U_2\right)
\left(\sup_{s<t,\, z\in U_2} p_E(s, z, y)\right)+ \E_x
\left[\tau_{U_1}\right] \left(\sup_{u\in U_1,\, z\in U_3}J(u,z)\right).
\end{equation}
\end{lemma}

\begin{lemma}\label{l:u_off1}
 There exists
 $c=c(\phi)>0$ such that  for every $(t,x,y) \in (0, \infty)
\times \bH \times \bH$,
$$
 p_{\bH}(t,x,y) \le c  \left(\sqrt{\frac{\Phi(\delta_{\bH}(x))}{t}} \wedge 1\right)
\left(\frac1{(\Phi^{-1}(t))^d}\wedge tJ(x, y)\right) .
$$
\end{lemma}

\pf
By \eqref{stssbound}, \eqref{e:J} and Lemma \ref{l:u_near}, it suffices to prove that
\begin{equation}
 p_{\bH}(t,x,y) \le c_1  \sqrt t \,
\sqrt{\Phi(\delta_{\bH}(x))}J(x, y) \quad \text{when } \delta_{\bH}(x) \le  \Phi^{-1}(t) \le |x-y|\label{e:dfnew}.
\end{equation}

We assume  $\delta_{\bH}(x) \le  \Phi^{-1}(t) \le |x-y|$ and let $x_0 =(\wt x, 0)$, $U_1:=B( x_0,
8^{-1}\Phi^{-1}(t)) \cap \bH$, $U_3:= \{z\in \bH: |z-x|>|x-y|/2\}$
and $U_2:=\bH\setminus (U_1\cup U_3)$.
Note that, by
Lemma \ref{l:phi-property} and Theorem \ref{L:2}(a), we have
\begin{equation}\label{e:upae1}
\E_{x}[\tau_{U_1}]\,\le\, \,
c_2\, \sqrt t \, \sqrt{\Phi(\delta_{\bH}(x))}.
\end{equation}
Since $U_1 \cap U_3 = \emptyset$ and
$
|z-x| > 2^{-1}|x-y|  \ge 2^{-1} \Phi^{-1}(t)$ for
$z\in U_3$,
we have for  $u\in U_1$ and $z\in U_3$,
\begin{eqnarray}\label{e:n001}
|u-z| \ge |z-x|-|x_0-x|-|x_0-u| \ge |z-x|- 4^{-1} \Phi^{-1}(t) \ge
\frac{1}{2}|z-x| \ge \frac{1}{4}|x-y|.
\end{eqnarray}
Thus, by \eqref{e:doubling-condition},
\begin{eqnarray}
\sup_{u\in U_1,\, z\in U_3}J(u,z) \le \sup_{(u,z):|u-z| \ge
\frac{1}{4}|x-y|}J(u,z)  \, \le \, c_3 j(|x-y|) .
\label{e:n01}
\end{eqnarray}
If $z \in U_2$,
\begin{equation}\label{e:one2}
\frac32 |x-y| \ge |x-y| +|x-z| \ge  |z-y| \ge |x-y| -|x-z| \ge
\frac{|x-y|}2 \ge 2^{-1} \Phi^{-1}(t).
\end{equation}
Thus, by \eqref{e:doubling-condition}, \eqref{stssbound} and \eqref{e:one2},
\begin{align}
&\sup_{s\le t,\, z\in U_2} p(s, z, y)
\le c_4 \sup_{ |x-y|/2 \leq |z-y|} t J(z,y)
\le c_5 t j(|x-y|) .\label{e:n02}
\end{align}
Applying Lemma \ref{l:gen1}, \eqref{e:upae1}, \eqref{e:n01}  and
\eqref{e:n02}, we obtain,
\begin{eqnarray*}
p_{\bH}(t, x, y) &\le& c_{6} \E_x[\tau_{U_1}]j(|x-y|) +
c_{6} \P_x\Big (X_{\tau_{U_1}}\in U_2 \Big)  t j(|x-y|) \\
&\le& c_{7} \sqrt t \, \sqrt{\Phi(\delta_{\bH}(x))} j(|x-y|) +
c_{6}\P_x\Big(X_{\tau_{U_1}}\in U_2\Big)  t j(|x-y|).
\end{eqnarray*}
Finally,
applying Lemma \ref{l2.1} and then \eqref{e:upae1}, we have
$$
\P_x \Big(X_{\tau_{U_1}}\in U_2 \Big) \le \P_x
\Big(X_{\tau_{U_1}}\in B( x_0, 8^{-1}\Phi^{-1}(t))^c\Big)
\,\le\, c_{8} \,
\frac1t\, \E_{x}[\tau_{U_1}]
\le c_9 t^{-1/2} \,
\sqrt{\Phi(\delta_{\bH}(x))} .
$$
Thus we  have proved
\eqref{e:dfnew}.
\qed

\begin{prop}\label{l:up3}
There exists
$c=c( \phi)>0$ such that for all $(t, x, y)\in (0,, \infty) \times
\bH\times \bH$,
\begin{eqnarray*}
 p_{\bH}(t, x, y)&\le& c \left(\sqrt{\frac{\Phi(\delta_{\bH}(x))}{t}} \wedge 1\right)\left(\sqrt{\frac{\Phi(\delta_{\bH}(y))}{t}} \wedge 1\right)
\left(\frac1{(\Phi^{-1}(t))^d}\wedge tJ(x, y)\right).
\end{eqnarray*}
\end{prop}

\pf
By Lemma \ref{l:u_off1} and
the lower bound of $p(t,x,y)$ in \eqref{stssbound},
there exists $c_1>0$ so that
 for every $z,w \in \bH$,  $
p_{\bH}(t/2, x, z)\le c_1(\sqrt{\Phi(\delta_{\bH}(x))/t} \wedge 1)
p(t/2, x,z).$
Thus, by the semigroup property and
the upper bound of $p(t,x,y)$ in \eqref{stssbound},
\begin{align*}
p_{\bH}(t,x,y)
&= \int_{\bH} p_{\bH}(t/2,x,z)p_{\bH}(t/2,z,y)dz \\
&\le c_1^2 \left(\sqrt{\frac{\Phi(\delta_{\bH}(x))}{t}} \wedge 1\right)\left(\sqrt{\frac{\Phi(\delta_{\bH}(y))}{t}} \wedge 1\right) \int_{\bH} p(t/2, x,z) p(t/2, y,z)dz\\
&\le c_1^2 \left(\sqrt{\frac{\Phi(\delta_{\bH}(x))}{t}} \wedge 1\right)\left(\sqrt{\frac{\Phi(\delta_{\bH}(y))}{t}} \wedge 1\right) p(t,x,y)\\
&\le c_2\left(\sqrt{\frac{\Phi(\delta_{\bH}(x))}{t}} \wedge 1\right)\left(\sqrt{\frac{\Phi(\delta_{\bH}(y))}{t}} \wedge 1\right)
\left(\frac1{(\Phi^{-1}(t))^d}\wedge tJ(x, y)\right).
\end{align*}
\qed

\begin{lemma}\label{L:4.2}
There exists $c=c(\phi)>0$ such that for any
$t >0$ and $y \in \R^d$,
$$
\P_y \left( \tau_{B(y, 8^{-1} \Phi^{-1}(t))} > t/3 \right)\,
\ge\, c.
$$
\end{lemma}

\pf By \cite[Proposition 4.9]{CK}, there exists $\eps=\eps(\phi)>0$ such that for every $t >0$,
$$
\inf_{y\in \bR^d } \P_y \left( \tau_{B(y, {16}^{-1} \Phi^{-1}(t)
)}
> \eps t \right) \ge \frac12.
$$
Suppose $\eps <\frac13 $, then by the parabolic Harnack inequality
in \cite{CK},
$$
c_1\,p _{B(y,8^{-1} \Phi^{-1}(t))}(\eps  t ,y,w) \, \le \,
p _{B(y,8^{-1} \Phi^{-1}(t))}(t/3 ,y,w)\qquad \hbox{for }
w \in B(y, {16}^{-1} \Phi^{-1}(t) ),
$$
where the constant $c_1=c_1(\phi)>0$ is
independent of $y\in \bR^d$. Thus
\begin{eqnarray*}
\P_y \left( \tau_{B(y,8^{-1} \Phi^{-1}(t))} > t/3\right) &=&
\int_{B(y, 8^{-1} \Phi^{-1}(t) )}
p_{B(y, 8^{-1} \Phi^{-1}(t) )}(t/3 ,y,w) dw\\
&\ge& c_1\int_{B(y, {16}^{-1} \Phi^{-1}(t) )} p_{B(y, 8^{-1} \Phi^{-1}(t)
)}(\eps t ,y,w) dw  \ge  \frac{c_1}2.
\end{eqnarray*}
\qed

The next result holds for any symmetric discontinuous
Hunt process that possesses a transition density and whose L\'evy system
admits a jumping density kernel. The proof is the same as that of \cite[Lemma 3.3]{CKS4} and so it is omitted here.

\begin{lemma}\label{l:gen}
Suppose that $U_1,U_2, U$ are open subsets of $\bR^d$ with $U_1,
U_2\subset U$ and ${\rm dist}(U_1,U_2)>0$. If $x\in U_1$ and $y \in
U_2$, then for all $t >0$,
\begin{equation}\label{eq:lb}
p_{U}(t, x, y)\,\ge\, t\,  \P_x(\tau_{U_1}>t)
\,\P_y(\tau_{U_2}>t)\inf_{u\in U_1,\, z\in U_2}J(u,z) \,.
 \end{equation}
\end{lemma}

\begin{lemma}\label{lower bound12}
 There exists $c=c(\phi)>0$ such that
for all $t >0$ and  $u, v\in \bR^d$ with $|u-v|\ge
\Phi^{-1}(t)/2$,
$$
 p_{B(u,\Phi^{-1}(t))\cup B(v,\Phi^{-1}(t))}(t/3, u, v)\,\ge\, c \,
t \, j(|u-v|).
$$
\end{lemma}

\pf Let $U= B(u,\Phi^{-1}(t))\cup B(v,\Phi^{-1}(t))$, $U_1=
B(u,\Phi^{-1}(t)/8)$,  $U_2=B(v,\Phi^{-1}(t)/8)$ and $K =\inf_{w\in
U_1,\, z\in U_2} j(|w-z|)$. We have by Lemma \ref{l:gen} that
\begin{eqnarray*}
p_{U}(t/3, u, v) \ge
 3^{-1} Kt \,  \P_u(\tau_{U_1}>t/3)\,
 \P_v(\tau_{U_2} >
t/3) \,.
\end{eqnarray*}
Moreover,
for $(w, z) \in U_1 \times U_2$,
$|w-z| \le |u-v| + |w-u|+ |z-v| \le |u-v| + \Phi^{-1}(t)/4
\le \frac32 |u-v|$.
Hence by \eqref{e:doubling-condition} $K \ge c_1j(|u-v|)$.
Thus
by Lemma \ref{L:4.2},
\begin{eqnarray*}
p_U(t/3, u,v)\ge    3^{-1} Kt  \,
 \left( \P_0(\tau_{B(0,\Phi^{-1}(t)/8)}>t/3) \right)^2
\ge
c_2 \,t \,j(|u-v|).
\end{eqnarray*}
\qed

\begin{lemma}\label{T:3.4}
Suppose that $D$ is an open subset of $\R^d$ and $(t, x, y)\in
(0, \infty) \times D\times D$ with $\delta_D(x) \ge \Phi^{-1}(t) \geq 2|x-y|$.
Then there exists  $c=c(\phi)>0$ such that
 \bee\label{e:lb1}
p_D(t,x,y) \,\ge\,c\, (\Phi^{-1}(t))^{-d}.
 \eee
\end{lemma}
\pf
Let  $t <\infty$ and $x, y \in D$ with $\delta_D(x) \ge
\Phi^{-1}(t) \geq 2|x-y|$.
By the parabolic Harnack inequality (\cite[Theorem 4.12]{CK}), there exists $c_1=c_1( \phi)>0$ such that
$$
p_D(t/2, x, w) \, \le  \, c_1\, p_D(t,x,y) \quad \mbox{for every
} w \in B(x, 2\Phi^{-1}(t)/3) .
$$
This together with Lemma \ref{L:4.2} yields that
\begin{eqnarray*}
p_D(t, x, y) &\geq & \frac{1}{c_1 \, | B(x,\Phi^{-1}(t)/2)|}
\int_{B(x,\Phi^{-1}(t)/2)} p_D(t/2, x, w)dw\\
&\geq &
c_2 (\Phi^{-1}(t))^{-d} \, \P_x \left( \tau_{B(x,\Phi^{-1}(t)/2)} >
t/2\right) \,\geq \, c_3 \, (\Phi^{-1}(t))^{-d},
\end{eqnarray*}
where $c_i=c_i(\phi)>0$ for $i=2, 3$.
\qed

For any $x\in \bH$ and $a, t>0$, we define
$
Q_x(a, t):=B((\wt x, 0), a\Phi^{-1}(t)) \cap \bH\ .
$

\begin{lemma}\label{l:keylon}
There exists $c=c (\phi)>0$ such that for all
$(t, x)\in (0, \infty) \times \bH$ with $\delta_{\bH}(x) < \Phi^{-1}(t)/2$,
$$
\P_x(\tau_{Q_x(2, t) }>t/3) \ge c
\frac{\sqrt{\Phi(\delta_{\bH}(x))}}{\sqrt t}.
$$
\end{lemma}

\pf
We fix $(t, x)\in (0, \infty) \times \bH$ with $\delta_{\bH}(x) < \Phi^{-1}(t)/2$.
The constants $c_1, \dots, c_8$ below are independent of  $t$ and $x$.
Without loss of generality we assume that $\wt x=\wt 0$ and
let $Q(a, t):=Q_0(a, t)$,
$x_1:=(\wt 0, \frac{3}{2}\Phi^{-1}(t))$ and $x_2:=(\wt 0, \frac{1}{4}\Phi^{-1}(t))$.
Note that, by L\'evy system and \eqref{e:doubling-condition},
\begin{align*}
&\P_{x_2}\Big(X_{\tau_{Q(1, t) }}\in B(x_1,  4^{-1}\Phi^{-1}(t))\Big)
\,\ge
\,
\P_{x_2}\Big(X_{\tau_{B(x_2,  4^{-1}\Phi^{-1}(t)) }}\in B(x_1,  4^{-1}\Phi^{-1}(t))\Big)\\
&=\int_{B(x_1,  4^{-1}\Phi^{-1}(t))} \int_{B(x_2,  4^{-1}\Phi^{-1}(t)) }G_{B(x_2,  4^{-1}\Phi^{-1}(t)) }(x_1, y) dy J(y,z)dz\\
&\ge c_1 \E_{0}[\tau_{B(0,  4^{-1}\Phi^{-1}(t))}]  \int_{B(x_1,  4^{-1}\Phi^{-1}(t))}  J(z)dz.
\end{align*}
Applying Theorem \ref{t:J-G}(a) and Lemmas \ref{l:phi-property} and \ref{l:exit-time}(b) to the above display, we get
$$
 \P_{x_2}\Big(X_{\tau_{Q(1, t) }}\in B(x_1,  4^{-1}\Phi^{-1}(t))\Big) \ge c_2
t\ |B(x_1,  4^{-1}\Phi^{-1}(t))|  \ t^{-1} \Phi^{-1}(t)^{d} \ge c_3\, .
$$
Thus, by Theorem \ref{L:2}(b),
\begin{eqnarray*}
\P_x\big(X_{\tau_{Q(1, t) }}\in B(x_1,  4^{-1}\Phi^{-1}(t))\big)
 \ge \ c_4
\P_{x_2}\big(X_{\tau_{Q(1, t) }}\in B(x_1,  4^{-1}\Phi^{-1}(t))\big)\tfrac{\sqrt{\Phi(\delta_{\bH}(x))}}{\sqrt{\Phi(\delta_{\bH}(x_2))}}
 \ge \ c_5
\tfrac{\sqrt{\Phi(\delta_{\bH}(x))}}{\sqrt t}.
\end{eqnarray*}
Now, using this,
Lemma \ref{L:4.2} and the strong Markov property,
\begin{eqnarray*}
&& \P_x\Big( \tau_{Q(2, t)} >t/3 \Big)\,\geq \, \P_x\Big( \tau_{Q(2, t)} >t/3, \
X_{\tau_{Q(1, t)}} \in B(x_1,  4^{-1}\Phi^{-1}(t)) \Big) \\
&& \ge \E_x\Big[ \P_{X_{\tau_{Q(1, t)
}}}\Big(\tau_{Q(2, t)} >t/3 \Big):
X_{\tau_{Q(1, t) }} \in B(x_1,  4^{-1}\Phi^{-1}(t))\Big]\\
&&\geq    \E_x\Big[ \P_{X_{\tau_{Q(1, t)
}}}
\Big(\tau_{B(X_{\tau_{Q(1, t)}}, \, 4^{-1}\Phi^{-1}(t))}>t/3 \Big)
: X_{\tau_{Q(1, t) }} \in B(x_1,  4^{-1}\Phi^{-1}(t))\Big]\\
&&=    \P_0
\Big(\tau_{B(0, \, 4^{-1}\Phi^{-1}(t))}>t/3 \Big)\P_x\Big(X_{\tau_{Q(1, t) }}\in B(x_1,  4^{-1}\Phi^{-1}(t))\Big)\\
&&\geq  c_7\P_x\Big(X_{\tau_{Q(1, t) }}\in B(x_1,  4^{-1}\Phi^{-1}(t))\Big) \ \ge \ c_8
\frac{\sqrt{\Phi(\delta_{\bH}(x))}}{\sqrt t}.
\end{eqnarray*}
This proves the lemma. \qed

Recall that $e_d$ denotes
the unit vector in the positive direction of the $x_d$-axis in $\R^d$. Now we are ready to prove the main result of this section
\begin{thm}\label{t:final}
There exists $c=c(\phi)>1$ such that for all
$(t, x,y)\in (0, \infty) \times \bH\times \bH$,
\begin{eqnarray*}
&&c^{-1} \left(\sqrt{\frac{\Phi(\delta_{\bH}(x))}{t}} \wedge 1\right) \left(\sqrt{\frac{\Phi(\delta_{\bH}(y))}{t}} \wedge 1\right)
\left(\frac1{(\Phi^{-1}(t))^d}\wedge tJ(x, y)\right)\\
&& \le p_{\bH}(t, x, y)\le c \left(\sqrt{\frac{\Phi(\delta_{\bH}(x))}{t}} \wedge 1\right) \left(\sqrt{\frac{\Phi(\delta_{\bH}(y))}{t}} \wedge 1\right)
\left(\frac1{(\Phi^{-1}(t))^d}\wedge tJ(x, y)\right) .
\end{eqnarray*}
\end{thm}

\pf
By Proposition \ref{l:up3}, we only need to show the lower bound of $p_{\bH}(t, x, y)$ in the theorem.
Fix $x,y\in \bH$.
Let $x_0=(\wt x, 0)$,  $y_0=(\wt y, 0)$,
$\xi_x:=x+32\Phi^{-1}(t)e_d$ and $\xi_y:=y+32\Phi^{-1}(t)e_d$.
If $\delta_{\bH}(x) < \Phi^{-1}(t)/2$,
by Lemmas \ref{L:4.2}, \ref{l:gen} and  \ref{l:keylon},
\begin{align*}
& \int_{B(\xi_x, 2\Phi^{-1}(t))} p_{{\bH}}(t/3,x,u)du\\
\ge& t\,  \P_x\left(
\tau_{Q_x(2, t)}>t/3\right)\left(\inf_{v\in Q_x(2, t) \atop
w\in B(\xi_x,4 \Phi^{-1}(t))}J(v,w) \right)
 \,\int_{B(\xi_x, 2\Phi^{-1}(t))}\P_u\left(\tau_{B(\xi_x, 4\Phi^{-1}(t))}>
t/3\right)du\\
\ge& c_1 t\,  \P_x\left(\tau_{Q_x(2, t)
}>t/3\right)
t^{-1}
(\Phi^{-1}(t))^{-d}
  \,   \P_0\left(\tau_{B(0, 8^{-1}\Phi^{-1}(t))}> t/3\right)
\,|B(\xi_x,2 \Phi^{-1}(t))| \nonumber\\
\ge& c_2 \P_x\left(\tau_{Q_x(2, t) }
>t/3\right)\ \ge \ c_3 \frac{\sqrt{\Phi(\delta_{\bH}(x))}}{\sqrt{t}}.
\end{align*}
On the other hand, if $ \delta_{\bH}(x) \ge \Phi^{-1}(t)/2$, by Lemmas \ref{L:4.2} and \ref{l:gen},
\begin{align*}
& \int_{B(\xi_x, 2\Phi^{-1}(t))} p_{{\bH}}(t/3,x,u)du  \nonumber
\\\ge& t\,  \P_x\left(\tau_{B(x,  8^{-1} \Phi^{-1}(t)) \cap \bH
}>t/3\right)\left(\inf_{v\in B(x_0, 2\Phi^{-1}(t)) \cap \bH \atop w\in
B(\xi_x,4 \Phi^{-1}(t))}J(v,w) \right) \,\int_{B(\xi_x,
2\Phi^{-1}(t))}\P_u\left(\tau_{B(\xi_x,
4\Phi^{-1}(t))}>t/3\right)du \nonumber\\
\ge& c_4 t\,  \P_x\left(\tau_{B(x,  8^{-1} \Phi^{-1}(t))}>t/3\right)
t^{-1}
(\Phi^{-1}(t))^{-d}\,  \P_0\left(\tau_{B(0, 8^{-1}\Phi^{-1}(t))}>
t/3\right) \,|B(\xi_x,2 \Phi^{-1}(t))| \nonumber\\
\ge& c_5
\P_0\left(\tau_{B(0,  8^{-1} \Phi^{-1}(t))}>t/3\right)^2  \ge c_6 .
\end{align*}
Thus
\begin{eqnarray}
\int_{B(\xi_x, 2 \Phi^{-1}(t))} p_{{\bH}}(t/3,x,u)du \ge c_7
\left(1\wedge \frac{\sqrt{\Phi(\delta_{\bH}(x))}}{\sqrt{t}} \right),
\label{e:loww_0}
\end{eqnarray}
and similarly,
\begin{eqnarray}
\int_{B(\xi_y, 2 \Phi^{-1}(t))} p_{{\bH}}(t/3,y,u)du \ge c_7
\left(1\wedge \frac{\sqrt{\Phi(\delta_{\bH}(y))}}{\sqrt{t}} \right).
\label{e:loww_01}
\end{eqnarray}

Now we deal with the cases $|x-y| \ge 5 \Phi^{-1}(t) $ and $|x-y| <
5\Phi^{-1}(t)$ separately.

\medskip

\noindent {\it Case 1}: Suppose that $|x-y| \ge  5 \Phi^{-1}(t) $.
Note that by the semigroup property  and Lemma~\ref{lower bound12},
\begin{align*}
&p_{\bH}(t,x,y)\nonumber\\
\geq& \int_{B(\xi_y, 2 \Phi^{-1}(t))}\int_{B(\xi_x, 2 \Phi^{-1}(t))}
p_{\bH}(t/3,x,u) p_{\bH}(t/3,u,v)p_{\bH}
(t/3,v,y)dudv \nonumber\\
\geq& \int_{B(\xi_y, 2 \Phi^{-1}(t))}\int_{B(\xi_x, 2 \Phi^{-1}(t))}
p_{{\bH}}(t/3,x,u)p_{B(u, \Phi^{-1}(t)) \cup
B(v,\Phi^{-1}(t))}(t/3,u,v)p_
{\bH}(t/3,v,y)dudv\nonumber\\
\geq& c_8 t \left(\inf_{(u,v) \in B(\xi_x, 2 \Phi^{-1}(t)) \times
B(\xi_y,2 \Phi^{-1}(t))} j(|u-v|)\right) \nn\\
&\quad \times\int_{B(\xi_y, 2
\Phi^{-1}(t))}
\int_{B(\xi_x, 2 \Phi^{-1}(t))}
p_{{\bH}}(t/3,x,u)p_{\bH}(t/3,v,y)dudv.
\end{align*}
It then follows from \eqref{e:loww_0}--\eqref{e:loww_01} that
\begin{align}\label{e:loww1}
p_{\bH}(t,x,
y)
\ge
 c_9  t \left(\inf_{(u,v) \in
B(\xi_x, 2 \Phi^{-1}(t)) \times B(\xi_y,2 \Phi^{-1}(t))} j(|u-v|) \right)
\left(\sqrt{\frac{\Phi(\delta_{\bH}(x))}{t}} \wedge 1\right)
\left(\sqrt{\frac{\Phi(\delta_{\bH}(y))}{t}} \wedge 1\right).
\end{align}
Using  the assumption $|x-y| \ge 5 \Phi^{-1}(t)$ we get that, for $u\in
B(\xi_x, 2 \Phi^{-1}(t))$ and $v\in B(\xi_y,2 \Phi^{-1}(t))$, $|u-v| \le
4 \Phi^{-1}(t) +|x-y| \le  2 |x-y|$. Hence
\begin{equation}\label{e:loww2}
\inf_{(u,v) \in B(\xi_x, 2 \Phi^{-1}(t)) \times B(\xi_y,2 \Phi^{-1}(t))}
j(|u-v|)\,\ge\, c_{10}
j(|x-y|).
\end{equation}
By \eqref{e:loww1} and \eqref{e:loww2}, we conclude that for $|x-y|
\ge 5  \Phi^{-1}(t)$
\begin{eqnarray*}
p_{\bH}(t, x, y)\ge c_{11} \left(\sqrt{\frac{\Phi(\delta_{\bH}(x))}{t}} \wedge 1\right)
\left(\sqrt{\frac{\Phi(\delta_{\bH}(y))}{t}} \wedge 1\right)t
j(|x-y|).
\end{eqnarray*}

\noindent {\it Case 2}: Suppose $|x-y| < 5 \Phi^{-1}(t)$. In this
case, for every $(u,v) \in  B(\xi_x, 2 \Phi^{-1}(t)) \times B(\xi_y,
2 \Phi^{-1}(t))$, $|u-v| \le 9 \Phi^{-1}(t)$. Thus, using the fact
that $\delta_{\bH}(\xi_x) \wedge \delta_{\bH}(\xi_y) \ge 32 \Phi^{-1}(t)$,
there exists $w_0 \in \bH$ such that
\begin{equation}\label{e:dfeggg}
B(\xi_x, 2 \Phi^{-1}(t)) \cup B(\xi_y, 2 \Phi^{-1}(t))
\subset B(w_0, 6 \Phi^{-1}(t))\subset B(w_0,
12 \Phi^{-1}(t))\subset \bH.
\end{equation}
Now, by the semigroup property and \eqref{e:dfeggg}, we get
\begin{align*}
&p_{\bH}(t,x,y)\nonumber\\
\geq& \int_{B(\xi_y, 2 \Phi^{-1}(t))}\int_{B(\xi_x, 2 \Phi^{-1}(t))}
p_{{\bH}}(t/3,x,u)p_{B(w_0,
12\Phi^{-1}(t))}(t/3,u,v)p_
{\bH}(t/3,v,y)dudv\nonumber\\
\geq&  \left(\inf_{u,v \in B(w_0, 6 \Phi^{-1}(t))} p_{B(w_0,
12\Phi^{-1}(t))}(t/3,u,v)\right) \int_{B(\xi_y, 2
\Phi^{-1}(t))}\int_{B(\xi_x, 2 \Phi^{-1}(t))}
p_{{\bH}}(t/3,x,u)p_{\bH}(t/3,v,y)dudv.
\end{align*}
It then follows from \eqref{e:loww_0}--\eqref{e:loww_01} and
Lemmas   \ref{l:phi-property} and \ref{T:3.4}
that
\begin{align*}
p_{\bH}(t, x, y) \ge c_{12} \left(\sqrt{\frac{\Phi(\delta_{\bH}(x))}{t}} \wedge 1\right)
\left(\sqrt{\frac{\Phi(\delta_{\bH}(y))}{t}} \wedge 1\right)(\Phi^{-1}(t))^{-d}
.
\end{align*}

Combining these two cases, we have proved the theorem.
\qed

Note that by using Theorem \ref{t:J-G}  we can express the sharp two-sided estimates for $p_{\bH}(t,x,y)$ solely in terms of the Laplace exponent $\phi$.

By integrating out time $t$ from the estimates in the preceding theorem, we can obtain
sharp two-sided  estimates of the Green function. Since the calculations are long and
somewhat cumbersome, we only state the result and omit the proof.
We refer the readers to \cite{CKS6} for  similar calculations.

\begin{thm}\label{t:final-green}
\begin{description}
\item{\rm (i)}
For all $d\ge 1$ there exists $c_1=c_1(d,\phi) >0$  such that for all
$(x, y)\in {\bH}\times {\bH}$,
$$
G_{\bH}(x, y)\ge c_1 \frac{\Phi(|x-y|)} {|x-y|^{d}}
 \left(1\wedge \frac{  \Phi (\delta_{\bH}(x))^{1/2}}{ \Phi(|x-y|)^{1/2}   }\right)
\left( 1\wedge \frac{ \Phi (\delta_{\bH}(y))^{1/2}}{ \Phi(|x-y|)^{1/2}} \right)
.
$$

\item{\rm (ii)}
If $d>(\delta_2 \vee \delta_4)$, then for all  $(x, y)\in
{\bH}\times {\bH}$,
$$
G_{\bH}(x, y)\asymp \frac{\Phi(|x-y|)} {|x-y|^{d}}
 \left(1\wedge \frac{  \Phi (\delta_{\bH}(x))^{1/2}}{ \Phi(|x-y|)^{1/2}   }\right)
\left( 1\wedge \frac{ \Phi (\delta_{\bH}(y))^{1/2}}{ \Phi(|x-y|)^{1/2}} \right)
.
$$

\item{\rm (iii)}
There exists
$c_2=c_2(d, \phi) >0$
such that for all
$(x, y)\in {\bH}\times {\bH}$ with $\Phi (\delta_{\bH}(x))  \Phi (\delta_{\bH}(y)) \le  \Phi(|x-y|)^2$,
$$
G_{\bH}(x, y)\le c_2
 \frac{ \Phi (\delta_{\bH}(x))^{1/2}  \Phi (\delta_{\bH}(y))^{1/2}}{ |x-y|^d}
.
$$

\item{\rm (iv)}
If $d=1$ and $\delta_1 \wedge \delta_3>1/2$, then  for all  $(x, y)\in
{\bH}\times {\bH}$,
$$
G_{\bH}(x, y)\asymp
\left( \frac{\Phi (\delta_{\bH}(x))^{1/2}  \Phi (\delta_{\bH}(y))^{1/2}}{\Phi^{-1}(\Phi (\delta_{\bH}(x))^{1/2}  \Phi (\delta_{\bH}(y))^{1/2})}\wedge \frac{ \Phi (\delta_{\bH}(x))^{1/2}  \Phi (\delta_{\bH}(y))^{1/2}}{ |x-y|} \right)
.
$$
\end{description}
\end{thm}


\medskip

{\bf Acknowledgements.}
We thank the referee for many helpful comments
on the first version of this paper.
We are especially grateful to the referee for suggesting
Lemma \ref{4.3.5} along with its proof which greatly simplifies the proof of Proposition \ref{l:main}.
\vspace{.1in}
\begin{singlespace}

\end{singlespace}
\end{doublespace}

\vskip 0.1truein

{\bf Panki Kim}

Department of Mathematical Sciences and Research Institute of Mathematics,

Seoul National University, Building 27, 1 Gwanak-ro, Gwanak-gu Seoul 151-747, Republic of Korea

E-mail: \texttt{pkim@snu.ac.kr}

\bigskip

{\bf Renming Song}

Department of Mathematics, University of Illinois, Urbana, IL 61801,
USA

E-mail: \texttt{rsong@math.uiuc.edu}

\bigskip

{\bf Zoran Vondra\v{c}ek}

Department of Mathematics, University of Zagreb, Zagreb, Croatia

Email: \texttt{vondra@math.hr}

\end{document}